\documentclass[12pt,a4paper]{article}
\usepackage{amsmath}
\usepackage{amssymb}
\usepackage{graphicx}
\usepackage[left=1in,right=1in,top=1in,bottom=1in]{geometry}
\usepackage{cite}
\usepackage[british]{babel}
 
\newcommand{\ud}{\mathrm{d}}
\newtheorem{theorem}{Theorem}[section]
\newtheorem{lemma}{Lemma}[section]

\newcommand{\rem}{\noindent \textbf{Remark. }}
\newcommand{\rems}{\noindent \textbf{Remarks. }}
\newcommand{\proof}{\noindent \textbf{Proof. }}
\newcommand{\qed}{$\square$}
\def\Exp{{\mathbb{E}}}
\def\Pr{{\mathbb{P}}}
\def\R{\mathbb{R}}
\def\Z{\mathbb{Z}}
\def\F{{\mathcal{F}}}

\def\L{{\mathcal{L}}}

\def\SS{\mathcal{S}}
\def\W{\mathcal{W}}
\def\N{\mathbb{N}}

\def\eps{{\varepsilon}}
\def\1{{\bf 1}}
\def\bx{{\bf x}}
\def\by{{\bf y}}
\def\bz{{\bf z}}
\def\be{{\bf e}}
\def\bo{{\bf o}}
\def\bq{{\bf q}}
\def\re{\mathrm{e}}
\def\0{{\bf 0}}
\def\H{{\mathcal{H}}}
\def\I{{\bf I}}
\def\M{{\bf M}}

\numberwithin{equation}{section}

\allowdisplaybreaks
 
\author{Iain M.\ MacPhee \qquad Mikhail V.\ Menshikov\\
{\footnotesize Department of Mathematical Sciences, Durham University,}\\
{\footnotesize South Road, Durham DH1 3LE, UK.}\\
\and Andrew R.\ Wade\\
{\footnotesize Department of Mathematics and Statistics, University of Strathclyde,}\\
{\footnotesize 26 Richmond Street, Glasgow G1 1XH, UK.}}

  \title{Moments of exit times from wedges for non-homogeneous
  random walks with asymptotically zero drifts}

\begin{document}

\maketitle

\begin{abstract}
We study quantitative asymptotics of planar random walks that are
spatially non-homogeneous but whose mean drifts have some regularity.
Specifically, we study the first exit time $\tau_\alpha$ from a
wedge with apex at the origin and interior half-angle $\alpha$
by a non-homogeneous 
random walk 
on $\Z^2$
with mean drift at $\bx$ of magnitude $O( \| \bx \|^{-1})$ as $\| \bx \| \to \infty$. 
This is the critical regime for the asymptotic behaviour:
under mild conditions, a previous result of the authors stated that $\tau_\alpha < \infty$ a.s.\ for any $\alpha$ (while for
a stronger drift field $\tau_\alpha$ is infinite with positive probability). 
Here we study the more difficult problem of the existence and non-existence
of 
moments $\Exp [ \tau_\alpha^s]$, $s>0$. Assuming (in common with much of the literature) a uniform bound on the walk's increments,
we show that for $\alpha < \pi/2$ there exists $s_0 \in (0,\infty)$ such that $\Exp [ \tau_\alpha^s]$ is finite for $s < s_0$ but infinite
for $s > s_0$; under specific assumptions on the drift field we show that  we can attain
$\Exp [ \tau_\alpha ^s] = \infty$ for any $s > 1/2$.
We show that for $\alpha \leq \pi$ there is a phase transition between  drifts of magnitude $O(\| \bx \|^{-1})$ (the {\em critical} regime)
and $o( \| \bx \|^{-1} )$ (the {\em subcritical} regime). In the subcritical regime
we obtain a non-homogeneous random walk
analogue of a theorem for Brownian motion due to Spitzer,
under considerably weaker conditions than those previously given (including work by Varopoulos) that assumed zero drift.
 \end{abstract}

\vskip 2mm


\vskip 2mm


\section{Introduction}
\label{intro}

By a random walk on $\R^d$ $(d \geq 2)$
we mean a discrete-time time-homogeneous
Markov
process on $\R^d$. If such a random walk is
spatially homogeneous, its position can be expressed as
a sum of i.i.d.~random vectors; such homogeneous random walks
are classical and have been extensively studied,
particularly when the state-space is $\Z^d$: see for example \cite{spitzerrw,lawlerbook}. 
The most subtle case is that of zero drift, i.e., when the increments
have mean zero. 

Spatial homogeneity, while simplifying the
mathematical analysis, is not always realistic for applications.
Thus it is desirable to study {\em non-homogeneous} random walks.
As soon as the spatial homogeneity assumption is relaxed, the
situation becomes much more complicated.
Even in the zero-drift case, a non-homogeneous random walk
can behave completely differently
to a zero-drift homogeneous random walk, and
can be transient in two dimensions, for instance.
This potentially wild behaviour  means that  techniques from the study of homogeneous random walks
 are difficult to apply. 

In this paper we continue the study of {\em angular asymptotics}, i.e.,
exit-from-cones problems, for non-homogeneous random walks that was started in \cite{mmw1}.
In \cite{mmw1} it was shown that, in contrast to recurrence/transience behaviour,
the angular properties of non-homogeneous random walks
are remarkably stable in some sense (as we describe later). 
We give more evidence to this effect in the present paper.

We study
non-homogeneous random walks with {\em asymptotically zero} mean-drift, that is,
the magnitude of the mean drift at $\bx \in \R^d$ tends to $0$ as $\| \bx \| \to \infty$.
This is the natural model in which to search for phase transitions
in asymptotic behaviour, as can be seen by analogy with the one-dimensional
problems considered by Lamperti \cite{lamp1,lamp2}, for instance.

Before formally defining our model and stating our theorems, we informally describe
existing results, the results in the present paper, and their significance.
In \cite{mmw1}, we studied the exit time $\tau_\alpha$ from a cone with interior half-angle $\alpha$ for a non-homogeneous random walk on $\Z^d$. 
For a zero-drift, {\em homogeneous} random walk, it is a classical result that $\tau_\alpha < \infty$ a.s.\ for any $\alpha$,
and tail asymptotics for $\tau_\alpha$ are known by comparison to a result of
Spitzer \cite{spitzer} for Brownian motion or by results of Varopoulos \cite{var1,var2}.
Our primary interest is how the situation changes when the walk is allowed to be non-homogeneous, and in particular, to quantify the effect
of introducing an asymptotically small  mean drift. 
 
 We will use $\mu (\bx)$ to denote the
one-step mean drift vector of the walk at $\bx$.
Unlike other asymptotic properties of random walk, it was shown
in \cite[Theorem 2.1]{mmw1} that the a.s.-finiteness of $\tau_\alpha$ remains valid
for non-homogeneous random walks {\em provided} $\| \mu (\bx) \| = O ( \| \bx \|^{-1} )$
as $\| \bx \| \to \infty$, under mild assumptions. In contrast, such a random walk can be
positive-recurrent, null-recurrent, or transient: see e.g.\ results of Lamperti \cite{lamp1,lamp2}.
On the other hand,
it was shown in \cite[Theorem 2.2]{mmw1} that a mean drift of magnitude
$\| \bx \|^{-\beta}$, $\beta \in (0,1)$, can ensure that the walk eventually remains in an arbitrarily
narrow cone: indeed, under mild conditions the walk is transient with a limiting direction
and a super-diffusive rate of escape \cite[\S 3.2]{mw3}.
 These facts motivate the following terminology. 
If $\| \mu  (\bx) \|$ is of magnitude 
(i) $o ( \| \bx \|^{-1} )$; (ii) $\| \bx \|^{-1}$;
(iii) $\| \bx \|^{-\beta}$, $\beta \in (0,1)$
we say that $\Xi$ is in the (i) {\em subcritical}; (ii) {\em critical};
(iii) {\em supercritical} regime, respectively.

The present paper is concerned with the critical and subcritical regimes.
Here we know from \cite{mmw1} that $\tau_\alpha < \infty$ a.s., but in the present
paper we are concerned with more detailed information about the random variable $\tau_\alpha$:
in particular, its {\em tails} (which moments do or do not exist). Thus the present paper
is concerned with {\em quantitative} information to complement the qualitative results of \cite{mmw1}.

There are two main themes of the present paper. 
First, we show that provided that $\| \mu  (\bx) \| = O (\| \bx \|^{-1} )$, $\tau_\alpha$ has a {\em polynomial}
tail, i.e., $\Exp[ \tau_\alpha^s]$ is finite for $s>0$ small enough but infinite for $s>0$ large enough.
Second, we demonstrate a {\em phase transition} in the tail behaviour of $\tau_\alpha$ between the
critical and subcritical regimes. 
Our main result on the subcritical regime will be that  not only does the
property  $\tau_\alpha < \infty$ a.s.~carry across
from the homogeneous zero-drift case, but also that 
finer information on the   moments 
of $\tau_\alpha$ also remains valid, under
mild additional conditions. On the other hand, we
give results that show that the
critical case is genuinely different: there is a quantitative phase transition
in the characteristics of $\tau_\alpha$
between the critical and sub-critical regimes.
 
Studying the moments of $\tau_\alpha$ is   much more difficult
  than determining whether $\tau_\alpha$ is a.s.\ finite, so in the present paper
we have to impose stronger conditions on the random walk than those in \cite{mmw1}. In particular,
to ease technical difficulties we impose a uniform bound on the increments of the walk (as opposed
to the 2nd moment bound used in \cite{mmw1}). The bounded increments assumption, although relatively strong, is prevalent
in the non-homogeneous random walk literature: see e.g.\ \cite{lawler91,mustapha,var1}. Moreover, we restrict
to two dimensions (in \cite{mmw1} the walk lived on $\Z^d$, $d \geq 2$). As well as again reducing technicalities, using 
 $\Z^2$
enables us to present our results as clearly as possible since even the Brownian motion case
becomes rather involved in higher dimensions \cite{burkh,dante}. We do not, however, need to assume
any symmetry for the increments (as required, for example, in \cite{lawler91,mustapha}).
 
 Before describing in detail our main results, we briefly survey some relevant literature.
In the homogeneous zero-drift setting,
for the analogous continuous problem
of planar Brownian motion in a wedge, a classical result of
Spitzer 
\cite[Theorem 2]{spitzer}
says that $\Exp [ \tau_\alpha^p] < \infty$ if and only if $p < \pi/(4 \alpha)$. A deep study of
  passage-time moments for Brownian motion in $\R^d$ was carried out by Burkholder \cite{burkh}.
The  random walk  problem has received 
less attention, even in the 
homogeneous zero-drift case. 
Varopoulos \cite{var1,var2} studied, using potential-theoretic methods, tails of passage-times for  zero-drift
random walks satisfying various conditions including bounded increments and isotropic covariance; 
some of the results of \cite{var1,var2} allow the walk to be
spatially inhomogeneous (at the expense of additional technical conditions, stronger than ours), but all require zero drift. 
From \cite{var1,var2} one can obtain a version of Spitzer's theorem for Brownian motion in the
case of zero-drift random walks satisfying appropriate regularity conditions. 
Exit times from cones for homogeneous random walks are also considered in \cite{garbit}. 
Other relevant results
specialize to the quarter-lattice $\Z^+ \times \Z^+$ 
\cite{cohen,klein} or the hitting-time of a half-line
\cite{lawlerbook,fukai}. Certain non-homogeneous random walks
with {\em linear} rate of escape were studied in \cite{flp}.
 
   A consequence of our results
 in the present paper is that Spitzer's theorem
 for Brownian motion
  essentially extends, under some moderate regularity conditions,
 to non-homogeneous  random walks with mean drifts that tend
 to zero   as $o(\|\bx\|^{-1})$. This considerably broadens the spectrum
 of random walks for which a Spitzer-type result is known; crucially, previous work
 has considered only the zero-drift case \cite{var1,var2}.
   
 We briefly comment on the techniques that we use
 in the present paper.
  Often it is possible to prove the existence
 of passage-time moments directly via semimartingale  (Lyapunov-type)
 criteria such as those in \cite{aim,lamp2}
 in the vein of Foster \cite{foster}. In the
 subcritical case for our non-homogeneous random walk, we have Lyapunov functions 
 that are well-adapted to do this. In the critical
 case, the non-homogeneity forces us to adopt a more direct
 approach, where nevertheless martingale ideas
 are central. The situation is similar for the 
 problem of non-existence of   moments,
 although even in the subcritical case rather delicate technical
 estimates are required. 
 
In the next section we formally define our model and state our main results.
 We also mention some possible directions for future research.

\section{Results and discussion}
\label{model}

We work in the plane $\R^2$; $\be_1, \be_2$
denote the standard orthonormal basis vectors and $\| \; \cdot \; \|$ the Euclidean norm. 
For $\bx \in \R^2$ we write
$\bx = (x_1,x_2)$ where $x_i = \bx \cdot \be_i$. Let $\0 = (0,0)$ denote
the origin. Our random walk will be
$\Xi = (\xi_t)_{t \in \Z^+}$, a Markov process whose state-space
is an unbounded subset $\SS$ of $\Z^2$. 

To
ensure that the walk cannot become trapped in lower-dimensional
subspaces  or finite sets,
we will assume the following  weak isotropy condition:

\begin{itemize}
\item[(A1)] There exist $\kappa>0$, $k  \in \N$ and
  $n_0 \in \N$ such that
\[ \min_{ \bx \in \SS; ~\by \in \{ \pm k  \be_1,\pm k  \be_2 \}} 
\Pr [ \xi_{t+n_0} - \xi_t = \by \mid \xi_t = \bx ] \geq \kappa ~~~(t \in \Z^+).\]
\end{itemize}

Note that (A1) is 
weaker than `uniform ellipticity' such as is often assumed
in the non-homogeneous random walk or random walk in random environment literature (see e.g.\
\cite{lawler91,mustapha}); for a discussion of the strength and implications of (A1), see \cite{mmw1}.
 
Let $\theta_t := \xi_{t+1} -\xi_t$ denote the jump of $\Xi$
at time $t \in \Z^+$. Since $\Xi$ is time-homogeneous and Markovian,
the distribution of the random vector $\theta_t$ depends only
upon the location $\xi_t \in \SS$ at time $t$. In other words,
 there exists a $\Z^2$-valued random field $\theta = (\theta(\bx))_{\bx \in \SS}$ 
such that for all $t \in \Z^+$,
\[  {\cal L} ( \xi_{t+1} - \xi_t \mid \xi_t  )  
= {\cal L} (\theta_t \mid \xi_t ) = {\cal L} ( \theta (\xi_t) ), \]
where $\L$ stands for `law'.
The law of $\theta$ is the {\em jump
distribution} of $\Xi$.
We write $\theta(\bx)$ in components as $(\theta_1(\bx),\theta_2(\bx))$.

Our second regularity condition is an assumption of
uniformly 
{\em bounded jumps}: 
\begin{itemize}
\item[(A2)]  There exists $b \in (0,\infty)$
such that
$\Pr[ \| \theta (\bx) \| > b ] = 0$ for all $\bx \in \SS$.
 \end{itemize}
It is likely that, as in \cite{mmw1}, this condition could be replaced
 by a moment assumption at the expense of
 some technical  work, but the assumption (A2) is frequently adopted in the
 literature: see e.g.\ \cite{mustapha,lawler91,var1}.

Under (A2), the moments of $\theta_t$ are well-defined.
Denote the one-step
{\em mean drift} vector
\[ \mu (\bx) := \Exp [ \theta_t \mid \xi_t = \bx ] = \Exp [ \theta (\bx) ] , \]
for $\bx \in \SS$, and write $\mu(\bx) = (\mu_1(\bx),\mu_2(\bx))$
in components. We are  
interested in the case of {\em asymptotically zero} mean drift, i.e.,
$\lim_{\| \bx \| \to \infty} \| \mu (\bx)\| =0$.

For $\alpha \in (0,\pi)$, we denote by $\W(\alpha)$ the (open) wedge
with apex at $\0$, principal axis in the $\be_1$ direction,
and interior half-angle $\alpha$:
\[ \W(\alpha ) := \{ \bx \in \R^d : \be_1 \cdot   \bx > \| \bx \| \cos \alpha \} .\]  
Thus $\W (\pi/4)= \{ (x_1,x_2) : x_1>0, |x_2| < x_1 \}$ is a quadrant
and $\W (\pi/2) = \{ (x_1,x_2) : x_1 >0 \}$ a half-plane.
The  case   $\alpha = \pi$ we will treat slightly differently:
for $s \geq 0$,
define 
\[  \H_s := \{ (x_1, x_2) : x_1 \leq 0, |x_2| \leq s \};\]
for $s>0$ this is a thickened half-line.
Then  with $b>0$ the jump bound in (A2),   set $\W (\pi) := \R^2 \setminus \H_b$.
(For convenience, we often call $\W(\pi)$ a `wedge' also.) It will
also be convenient to
set $\W (\alpha) := \R^2$ for any $\alpha > \pi$.

Our primary quantity is
the random walk's first exit time from the wedge $\W(\alpha)$. 
With the usual convention that $\min \emptyset := \infty$,
define 
\begin{align}
\label{tau} \tau_\alpha   
 := \min \{ t \in \Z^+ : \xi_t \notin  \W (\alpha )
    \}.\end{align}
    
The following fundamental result says that as soon as the mean drift decays fast enough,
$\tau_\alpha$ is a.s.\ finite. Theorem \ref{thm1} is essentially contained in \cite{mmw1}:
indeed, \cite[Theorem 2.1]{mmw1} gives such a result under conditions much weaker than (A2) and
in general dimensions $d \geq 2$, but not including the case $\alpha = \pi$ (hitting the thickened
half-line). We will give a self-contained proof of Theorem \ref{thm1}
that requires minimal extra work on top of that  to obtain the main results
of the present paper.
 
\begin{theorem}
\label{thm1}
Suppose that (A1) and (A2)
 hold, and
that for $\bx \in \SS$ as $\|\bx\| \to \infty$,
\begin{align}
\label{drift1}
\| \mu (\bx)\| = O ( \| \bx \|^{-1} ) .\end{align}
Then for any $\alpha \in (0,\pi]$ 
 and any $\bx \in \W (\alpha)$, $\Pr [ \tau_\alpha < \infty \mid  \xi_0 = \bx ]  = 1$.
 \end{theorem}

Our first substantially new result, Theorem \ref{thm2},
 gives information on the tails of  $\tau_\alpha$,
 $\alpha < \pi/2$.
 In particular, it shows that even for this non-homogeneous
 walk, the tail behaviour
 is essentially polynomial
 in character, as in the zero-drift case: compare Theorem \ref{thm5} below. However,
 the `heaviness' of the tail (i.e., the exponent $s_0$ in the statement
 of Theorem \ref{thm2}) will depend on the details of the walk: compare
 Theorems \ref{thm9} and \ref{thm5} below. For a one-dimensional
 analogue of this result, see the Appendix in \cite{aim}.
 
 \begin{theorem}
 \label{thm2}
 Suppose that (A1) and (A2) hold,
 $\alpha \in (0,\pi/2)$,
 and that  for $\bx \in \SS$,
 (\ref{drift1}) holds  as $\|\bx\| \to \infty$. Then  there exist $s_0, A \in (0,\infty)$
 such that:
 \begin{itemize}
 \item[(i)] if $s < s_0$, then  $\Exp [ \tau_\alpha^s \mid \xi_0 = \bx ] < \infty$ for any $\bx \in \W(\alpha)$;  
 \item[(ii)] if $s > s_0$, then   $\Exp [ \tau_\alpha^s \mid \xi_0 = \bx ] = \infty$ for any $\bx \in \W(\alpha)$
 with $\| \bx \| \geq A$.
 \end{itemize}
 \end{theorem}
\rems
(a) It is an open problem to show that  
Theorem \ref{thm2} holds for $\alpha \geq \pi/2$.\\
\noindent
(b) Theorem \ref{thm2}(ii) cannot be strengthened
to {\em all} $\bx \in \W(\alpha)$ without  stronger regularity conditions
on the walk $\Xi$. Indeed, under (A1), it may be that for $\xi_t$ close to the boundary
of $\W(\alpha)$, $\xi_{t+1}$ is outside $\W(\alpha)$ with probability 1; however, this cannot
occur for $\|\xi_t\|$ large enough by our asymptotically zero drift assumption: see Lemma \ref{lem88} below.
The same remark applies to  our other non-existence of moments results that follow.\\

Walks satisfying
Theorem \ref{thm2}
can have radically different characteristics.
For example, for small enough wedges 
a zero-drift walk will have $\Exp [ \tau_\alpha ] <\infty$
(see Theorem \ref{thm5} below). On the other
hand, the next result implies that for any $\alpha \in (0,\pi/2)$,
for
a suitably strong $O(\| \bx\|^{-1})$
drift field,  $\Exp [ \tau_\alpha ] =\infty$.
In fact, Theorem \ref{thm9} says that for any $\eps>0$,
there exist walks satisfying the conditions
of Theorem \ref{thm2} for which  $(1/2)+\eps$
moments of $\tau_\alpha$ do not exit.
An open question is to determine
whether $(1/2)-\eps$ moments
can be infinite
under the conditions of Theorem \ref{thm2}.

We take the random walk
to have dominant drift
in the principal direction. Specifically,
we assume that
there exists $c >0$
for which
\begin{align}
\label{qq2}
\liminf_{\| \bx \| \to \infty} ( \| \bx \| \mu_1 (\bx) ) \geq c, 
~~~
\lim_{ \| \bx \| \to \infty } ( \| \bx \| \mu_2 (\bx) ) = 0.
 \end{align} 

\begin{theorem}
\label{thm9}
Suppose that (A1) and (A2) hold, and
$\alpha \in (0,\pi/2)$.
Then for any $s>0$, there exist $c_0, A \in (0,\infty)$
such that if (\ref{qq2}) holds for any $c > c_0$, then
for all $\bx \in \W(\alpha)$ with $\| \bx \| \geq A$, 
$\Exp [ \tau_\alpha^{(1/2)+s} \mid \xi_0 = \bx ] = \infty$.
\end{theorem}

Our final result, Theorem \ref{thm5}, gives
sharp tail asymptotics
for $\tau_\alpha$
in the {\em subcritical}
regime. To obtain such a sharp result, we need
to assume additional regularity for $\Xi$: specifically, we need
to control the covariance structure of the increments of the walk.
Denote
the covariance matrices
 $\M=(M_{ij})_{i,j \in \{1,2\}}$ of $\theta$ by
\[ \M (\bx) := \Exp [ {\theta_t}^{\!\!\top} \theta_t \mid \xi_t = \bx]
= \Exp [\theta (\bx)^\top \theta (\bx) ] ,\]
for $\bx \in \SS$, where $\theta_t$ is viewed as  a row-vector.
When (A1) holds, $\Pr [ \xi_{t+1} \neq \bx \mid \xi_t = \bx ]$ is uniformly positive 
\cite[p.\ 4]{mmw1} so that
$M_{11}(\bx)+M_{22}(\bx) >0$ uniformly in $\bx$.

Theorem \ref{thm5}  
shows that the critical
exponent 
for the moment problem
depends only on $\alpha$ and 
is the same in
this random walk setting
as in the Brownian motion
case, where the result
is due to Spitzer \cite[Theorem 2]{spitzer}. In 
particular,
Theorem \ref{thm5}
includes the case  
of a homogeneous random walk with zero  drift,
 where the result follows from  \cite[Theorem 4]{var1} (see also \cite{var2}).
We write $\bo(1)$ for a $2 \times 2$
matrix each of whose entries
is $o(1)$. 

\begin{theorem}
\label{thm5}
Suppose that (A1) and (A2) hold,
and
there exists $\sigma^2 \in (0,\infty)$
such that 
\begin{equation}
\label{sub2}
 \|\mu(\bx)\| = o( \| \bx \|^{-1} ), \textrm{ and }  
\M (\bx) =  \sigma^2  \I + \bo (1),\end{equation}
as $\| \bx \| \to \infty$. Suppose that $\alpha \in (0,\pi]$.
\begin{itemize}
\item[(i)] If $s \in [0, \pi/(4\alpha))$
and  $\bx \in \W (\alpha)$, $\Exp [ \tau^s_\alpha \mid \xi_0 = \bx] < \infty$.
\item[(ii)]
If
$s> \pi/(4\alpha)$
and  $\bx \in \W (\alpha)$ with $\| \bx \|$ sufficiently large,
$\Exp [ \tau^s_\alpha \mid \xi_0 = \bx] = \infty$.
\end{itemize}
\end{theorem}

Certain cases
of Theorem \ref{thm5}
extend results of
Klein Haneveld and Pittenger \cite{klein}
and Lawler \cite{lawlerbook}
for homogeneous zero-drift random
walks (i.e., sums
of i.i.d.~mean-zero random vectors)
to non-homogeneous random
walks with small drifts. First,
for hitting
a half-line
($\alpha =\pi$),
Theorem \ref{thm5} 
implies that 
$1/4$-moments
are critical,
a result
obtained
for homogeneous
zero-drift 
random walks by Lawler
(see (2.35) in \cite{lawlerbook}, also \cite{fukai}).
 Second,
in the case
of a quadrant ($\alpha =\pi/4$),
Theorem \ref{thm5}(ii)
implies that $\Exp [\tau_{\pi/4}^s] =\infty$
for $s>1$,
a result contained
in   \cite[Theorem 1.1]{klein}
for a homogeneous
zero-drift random walk
with certain regularity conditions (see also \cite[Theorem 1.1]{cohen}).
 
 The outline of the rest of the paper is as follows.  Section \ref{prelim} collects some preparatory
results. Section \ref{prf1} is devoted to the critical case  and   the
proofs of Theorems \ref{thm1}, \ref{thm2} and \ref{thm9}, while Section \ref{prf2}
is devoted to the subcritical case and the proof of Theorem \ref{thm5}. The proofs
of the existence and non-existence of moments results are largely separate.

 \section{Preliminaries}
 \label{prelim}
  
\subsection{Semimartingale criteria}
\label{semi}

In this section we collect some general semimartingale-type results that we
need. 
Let $(\F_t)_{t \in \Z^+}$
be a filtration on a probability
space $(\Omega,\F,\Pr)$.
Let $(Y_t)_{t \in \Z^+}$
be a discrete-time $(\F_t)_{t \in \Z^+}$-adapted
stochastic process
taking values in  $[0,\infty)$. Typically,
when we come to apply the following lemmas
later on, we will have $Y_t = r (\xi_t)$
for some $r: \R^2 \to [0,\infty)$.

Following
work of Lamperti \cite{lamp2},
the primary results
available for establishing the existence and non-existence
of passage-time moments for a (not necessarily Markov) stochastic
process are contained in \cite{aim}. For some of the applications
in the present paper, we could not apply these general
results and so have to use other techniques.
 
The following existence result is contained
in Theorem 1 of \cite{aim}. 
 
\begin{lemma} \label{aimthm}
Let $(Y_t)_{t \in \Z^+}$ be an $(\F_t)_{t\in\Z^+}$-adapted
stochastic process taking values in an unbounded
subset
  of $[0,\infty)$. For
 $B>0$ set $\upsilon_B := \min \{ t \in \N : Y_t \leq B \}$. Suppose that
there exist $C, p_0 \in (0,\infty)$ such that for any $t \in \Z^+$,
 $Y_t^{2p_0}$
is integrable, and
\[ \Exp [ Y_{t+1}^{2p_0} - Y_t ^{2p_0} \mid \F_t ] \leq 
-C
  Y_t ^{2p_0-2} , ~\textrm{on}~ \{ \upsilon_B > t \} .\]
Then for any $p \in [0,p_0)$, 
 for any $x$,
$\Exp [ \upsilon_B^p \mid Y_0 = x]   
< \infty$.
 \end{lemma}

The corresponding non-existence
result that we will need is
 Corollary 1 in \cite{aim}:

\begin{lemma} \label{aimthm2}
With the notation
of Lemma \ref{aimthm}, 
suppose that
there exist $C, D , p_0 \in (0,\infty)$ and $r>1$ such that for any $t \in \Z^+$
the following 3 conditions hold on $\{ \upsilon_B > t\}$:
  \begin{align}
  \label{noncon1}
 \Exp [ Y_{t+1}^{2p_0} - Y_t ^{2p_0} \mid \F_t ] & \geq 
0; \\
\label{noncon2}
\Exp[ Y_{t+1}^2 - Y_t^2 \mid \F_t ] & \geq  - C; \\
\label{noncon3}
\Exp [ Y_{t+1}^{2r} - Y_t^{2r} \mid \F_t ] & \leq D Y_t^{2r-2} .\end{align}
Then 
 for any $p > p_0$, 
for any $x$ large enough, 
$\Exp [ \upsilon_B^p \mid Y_0 = x ] = \infty$.
\end{lemma}

\subsection{Lyapunov functions}
\label{harm}

In this section we introduce some Lyapunov  functions that we will use
to study our random walk  in the subcritical case,
and analyze their basic properties. These functions will be built upon standard harmonic
functions in the plane, as were employed by Burkholder \cite{burkh}
in his sharp analysis of the exit-from-cones problem for  Brownian motion;
it is natural that they are the correct tools when our random walk
is sufficiently close to zero-drift. We need some more notation.

For $\bx =(x_1,x_2) \in \R^2$ we use polar coordinates $(r,\varphi)$
relative to the ray $\Gamma_0$
in the $\be_1$ direction
 starting at $\0$.
Thus if $r = \| \bx \|$ and $\varphi \in (-\pi,\pi]$ is the angle, measuring anticlockwise, of the ray through
$\0$ and $\bx=(x_1,x_2)$ from the ray $\Gamma_0$,
we have $x_1 = r \cos \varphi$ and
$x_2= r \sin \varphi$. 
We occasionally write $\varphi$ as $\varphi(\bx)$ for clarity.
Let $\be_r (\varphi) = \be_1 \cos \varphi
+ \be_2 \sin \varphi$, the radial unit vector,
and
$\be_\perp (\varphi) = -\be_1 \sin \varphi + \be_2 \cos \varphi$,
the transverse unit vector. Note that in polar coordinates,
$\W(\alpha)  = \{ \bx \in \R^2 : r>0, -\alpha < \varphi < \alpha \}$.

Let $B_r(\bx)$ denote the closed
Euclidean  ball (a disk) of radius $r$
centred at $\bx \in \R^2$. For $\alpha \in (0,\pi]$
and  $s \geq 0$ define
the modified wedge 
\[ \W_s (\alpha) := \W (\alpha) \setminus B_s (\0) = \{ \bx \in \W (\alpha) : \| \bx \| > s\},\]
which is $\W (\alpha)$ with a disk-segment around the origin 
removed.
During our proofs, 
we will often
 work with
the exit time from the locally modified set $\W_A(\alpha)$ for some fixed
(large) value of $A>0$.
Let $\W_0(\alpha):=\W(\alpha)$ and
for $A \geq 0$, define  
\begin{align}
\label{taua}
  \tau_{\alpha,A} 
 := \min \{ t \in \Z^+ : \xi_t \notin  \W_A (\alpha )
     \}.\end{align}
     Then 
     $\tau_{\alpha,A} \geq \tau_{\alpha,B}$ for $B \geq A$, and
     $\tau_{\alpha,0} = \tau_\alpha$ with the notation
     of (\ref{tau}).

We will use multi-index notation for partial derivatives on $\R^2$. For 
$\sigma = (\sigma_1,\sigma_2) \in \Z^+ \times \Z^+$, $D_\sigma = D_{\sigma_1 \sigma_2}$ will 
denote $D_1^{\sigma_1} D_2^{\sigma_2}$
where $D_j^k$ for $k \in \N$ is $k$-fold 
differentiation with respect to $x_j$, and $D^0_j$ is the identity
operator. We also use the notation $|\sigma|:=\sigma_1+\sigma_2$
and $\bx^\sigma :=  x_1^{\sigma_1} x_2^{\sigma_2}$.
 
For $w>0$, define the function $f_w :\R^2 \to \R$ by
\begin{align}
\label{fdef} f_w (\bx) := f_w (r,\varphi) := r^w \cos ( w \varphi ).\end{align}
Differentiating, using the appropriate from of the chain rule, shows that for any $w>0$,
\begin{align}
\label{diff1}
D_1 f_w (r,\varphi) = w r^{w-1} \cos ( (w-1) \varphi) ; ~~~
D_2 f_w (r,\varphi) = -w r^{w-1} \sin ( (w-1) \varphi) ,\\
\label{diff2}
\text{and} ~~~ D_1 D_2 f_w (r,\varphi)
= D_2 D_1 f_w (r,\varphi)
= w (w-1) r^{w-2} \sin ((w-2) \varphi).\end{align}
Moreover,  
$f_w$ is   harmonic on $\R^2$, since
\begin{align}
\label{diff3}
D_1^2 f_w (r,\varphi) = w(w-1) r^{w-2} \cos ((w-2) \varphi) = - D_2^2 f_w (r,\varphi).\end{align}
For $w>1/2$, $f_w$ is positive in the interior of the
wedge $\W (\pi/(2w))$, and $0$ on the boundary $\partial \W(\pi/(2w))$;
$f_{1/2}$ is positive on $\R^2 \setminus \H_0$ and zero on
the  half-line $\H_0$. For $w \in (0,1/2)$, $f_w$ is positive
throughout $\R^2$.
As an example,  the harmonic function
\begin{align}
\label{qq1}
f_2 (\bx) = r^2 \cos (2 \varphi) = x_1^2-x_2^2 \end{align}
 is positive on the quadrant $\W(\pi/4)$
and zero on $\partial \W(\pi/4)$. 

It follows by repeated applications of the chain rule
that $f_w$ and all of its derivatives $D_\sigma f_w$ are of the 
form $r^k u (\varphi)$ where $u$ is bounded, and hence for any
$\sigma$ with $|\sigma|=j$ there exists a constant $C  \in (0,\infty)$ such that
for all $\bx \in \R^2$,
\begin{align}
\label{az1}
- C r^{w-j} < D_\sigma f_w (\bx) < C r^{w-j} .\end{align} 

The next result gives expressions
for the first three moments of the
jumps of $f_w (\xi_t)$.

\begin{lemma}
\label{fmoms}
Suppose that (A2) holds.
Then with $f_w$ defined at (\ref{fdef}), for $w>0$, 
 there exists $C \in (0,\infty)$ such that
  for any $\bx \in \SS$,  
 \begin{equation}
 \label{fbnd}
 \Pr [ | f_w (\xi_{t+1}) - f_w (\xi_t) |
  \leq C (1+ \| \bx \|)^{w-1} \mid \xi_t = \bx ]
 =1 .\end{equation}
 Also,
 for any $\bx \in \SS$
 as $r= \| \bx \| \to \infty$,
we have the following asymptotic expansions:
\begin{align}
\label{fmom1}
 \Exp [ f_w(\xi_{t+1}) - f_w(\xi_t) \mid \xi_t = \bx ] 
  = {} & w r^{w-1} \left( \mu_1(\bx) \cos  ( (w-1) \varphi ) 
 -\mu_2(\bx) \sin ((w-1)\varphi) \right)  
  \nonumber\\ 
&
{}+ \frac{1}{2} \left(M_{11} (\bx) -M_{22} (\bx) \right) w(w-1) r^{w-2} \cos((w-2)\varphi) \nonumber\\
& {}+
M_{12} (\bx) w (w-1) r^{w-2} \sin ((w-2)\varphi) 
+  O(r^{w-3} );\\
\label{fmom2}
   \Exp [ ( f_w(\xi_{t+1}) - f_w(\xi_t))^2 \mid \xi_t = \bx ] 
  = {} & w^2    r^{2w-2} \left( M_{11}(\bx) \cos^2 ((w-1)\varphi) + M_{22}(\bx) \sin^2 ((w-1)\varphi) 
      \right) \nonumber\\
      & {}-   M_{12} (\bx) w^2 r^{2w-2} \sin ( 2(w-1) \varphi) + O(r^{2w-3}) ;\\
\label{fmom3}
 \Exp [ ( f_w(\xi_{t+1}) - f_w(\xi_t))^3 \mid \xi_t = \bx ] 
   = {} & O( r^{3w-3}).\end{align}
\end{lemma}
\proof Since $f_w$ is smooth, Taylor's theorem
with Cartesian coordinates and the Lagrange form for the remainder applied in a disk of radius $b$ at any $\bx \in \R^2$ implies
\begin{align*}
f_w (\bx + \by) = f_w (\bx) + \sum_j y_j (D_j f_w ) (\bx + \eta \by),
\end{align*}
for some $\eta = \eta (\by) \in (0,1)$,
for any $\by =(y_1, y_2)$ with $\| \by \| \leq b$.
Taking $\by = \xi_{t+1} - \xi_t = \theta_t$,
conditioning on $\xi_t =\bx$,
we then obtain, with (\ref{diff1}), a.s.,
for some $C \in (0,\infty)$,
\[ | f_w (\xi_{t+1}) - f_w(\xi_t) | \leq C \| \bx + \eta \theta (\bx) \|^{w-1} ,\]
for any $\by \in \Z^2$.
Now (A2) implies (\ref{fbnd}).
 
For the moment estimates, we 
include more terms in the
Taylor expansion to obtain
\begin{align}
\label{jjj}
& \Exp [ f_w (\xi_{t+1}) - f_w (\xi_t) \mid \xi_t = \bx ] = \sum_j \Exp [ \theta_j ( \bx) ] (D_j f_w) (\bx) 
+ \frac{1}{2} \sum_j  \Exp [\theta_j (\bx) ^2  ] (D_j^2 f_w) (\bx) \nonumber\\
  & \qquad {}+\sum_{i < j} \Exp [ \theta_i (\bx) \theta_j (\bx) ]  (D_i D_j f_w)
(\bx ) 
+ \frac{1}{6} \Exp \left[   \sum_{\sigma : |\sigma| =3} \theta^\sigma
(\bx) (D_\sigma f_w ) (\bx + \eta 
\theta (\bx) )    \right] ,\end{align}
for some $\eta = \eta (\theta_t) \in (0,1)$. 
By (A2),
 $\Exp [\theta_j (\bx) ^3
 ] = O(1)$,
 so that
using (\ref{az1}) 
the final term
in (\ref{jjj}) is $O(r^{w-3})$. Then using
(\ref{diff1}), (\ref{diff2}) and (\ref{diff3})
in (\ref{jjj}), we obtain (\ref{fmom1}).
 
In a similar fashion, we obtain (\ref{fmom2}). Specifically,
\begin{align*}
\Exp [ ( f_w (\xi_{t+1}) - f_w (\xi_t))^2 \mid \xi_t = \bx ] 
= \sum_j  \Exp [\theta_j (\bx)^2
 ] (( D_j f_w ) (\bx))^2
 \\
{}+ 2 \sum_{i<j}  \Exp [ \theta_i (\bx) \theta_j (\bx) ] 
 (D_i f_w)(\bx)(D_j f_w) (\bx)  + O( r^{2w-3} ) ,\end{align*}
using (A2), and (\ref{fmom2}) follows using (\ref{diff1}).
Finally, 
\begin{align*}
\Exp [ ( f_w (\xi_{t+1}) - f_w (\xi_t))^3 \mid \xi_t = \bx ]
 = \sum_j \Exp [\theta_j (\bx) ^3
]
(( D_j f_w ) (\bx))^3
  + O(r^{3w-4})
,\end{align*}
and by (A2),
$\Exp [\theta_j (\bx) ^3]
  = O(1)$. Then (\ref{fmom3}) follows from (\ref{diff1}). \qed\\
  
When $\Xi$ has zero drift, one expects that $(f_w(\xi_t))_{t \in \Z^+}$ is `almost' a martingale,
 keeping the Brownian  analogy in mind  \cite{burkh}. 
Thus the process $(f_w(\xi_t))_{t \in \Z^+}$ will be useful  when $\| \mu (\bx) \| = o(\| \bx \|^{-1})$. 
In order to apply the semimartingale
criteria of Section \ref{semi}, 
we often 
want to modify our process $(f_w(\xi_t))_{t \in\Z^+}$, 
to obtain either a submartingale or a supermartingale. So,
in Lemma \ref{fgmoms} below,
 we study the process
$(f_w(\xi_t)^{\gamma})_{t \in \Z^+}$ where $\gamma \in \R$.
Recall that 
for $w <   \pi/(2\alpha)$, $f_w (\bx)$ is positive on a wedge $\W (\pi/(2w))$
 bigger than $\W (\alpha)$. The following result is simple but important.
 \begin{lemma}
   \label{lbdlem}
 Suppose that $\alpha \in (0,\pi]$ and  $w \in (0,\pi/(2\alpha))$.
 Then there exists $\eps_{\alpha,w} = \cos (w \alpha ) >0$ such that
 for all $\bx \in \W(\alpha)$,
  \begin{align}
  \label{lbd}
   \eps_{\alpha, w}  r^w \leq  f_w (\bx) \leq   r^w .\end{align}
  Moreover, for $k \geq 0$ we have that if $w \geq k/2$
  then for all $\bx \in \W(\alpha)$, 
  \begin{align}
  \label{lbd2}
   \cos ((w-k) \varphi) \geq \eps_{\alpha,w} >0  .\end{align}
     \end{lemma}
  \proof 
  For fixed
  $\alpha\in (0,\pi]$ and fixed $w \in (0,\pi/(2\alpha))$ we have
  \[ \eps_{\alpha,w} := \inf_{\bx \in \W(\alpha)} \cos ( w \varphi)
  = \inf_{\varphi \in (-\alpha,\alpha)} \cos (w \varphi) = \cos ( w \alpha ) > 0 ,\]
  since $w\alpha \in (0,\pi/2)$. Then   (\ref{lbd}) follows from (\ref{fdef}).
   The statement (\ref{lbd2})
  follows similarly, using  the fact that for $w \geq k/2$ and $k \geq 0$,
$-w \leq -k/2 \leq w - k \leq w$, 
  so  
  \[ \inf_{\bx \in \W(\alpha)} \cos ((w-k) \varphi) 
  \geq \inf_{\varphi \in (-\alpha,\alpha)}
  \cos ( w \varphi ) = \eps_{\alpha,w} .\]
  This completes the proof. \qed

\begin{lemma}
\label{fgmoms}
Suppose that (A2) holds.
Suppose that $\alpha \in (0,\pi]$, $\gamma \in \R$, and $w \in (0,\pi/(2\alpha))$.
Then for all
$\bx \in \W (\alpha)$  we have
that as $r = \| \bx \| \to \infty$,
\begin{align*}
& \Exp [ f_w(\xi_{t+1})^{\gamma} - f_w(\xi_t)^{\gamma} \mid \xi_t = \bx ] \nonumber\\
= &~ 
\gamma f_w(\bx)^{\gamma-1}  w r^{w-1} \left(
\mu_1(\bx) \cos((w-1)\varphi) - \mu_2(\bx) \sin((w-1)\varphi) \right)   \nonumber\\
& {}+ \gamma f_w (\bx)^{\gamma-1} M_{12} (\bx) w(w-1) r^{w-2} \sin((w-2)\varphi)
   \nonumber\\
&  {}+ \frac{1}{2} \gamma f_w (\bx)^{\gamma-1} \left( M_{11}(\bx)-M_{22}(\bx) \right)
  w(w-1) r^{w-2} \cos ((w-2)\varphi) \nonumber\\
&  {}+ \frac{1}{2} \gamma (\gamma-1) f_w (\bx)^{\gamma-2} w^2 r^{2w-2}
  \left( M_{11} (\bx)\cos^2 ((w-1)\varphi)
  +M_{22}(\bx) \sin^2 ((w-1)\varphi) \right) \nonumber\\
&  {}- \frac{1}{2} \gamma (\gamma-1) f_w(\bx)^{\gamma-2}
  w^2 r^{2w-2} M_{12} (\bx) \sin(2(w-1)\varphi) 
  + O(f_w (\bx)^{\gamma-3} r^{3w-3} )  .\end{align*}
\end{lemma}
\proof
Let $\Delta := f_w(\xi_{t+1}) - f_w(\xi_t)$. 
Then for $\gamma \in \R$ and $\bx \in \W(\alpha)$,
\[ \Exp [ f_w(\xi_{t+1})^\gamma - f_w(\xi_t)^\gamma \mid \xi_t = \bx]
= f_w(\bx)^\gamma \Exp \left[ \left.
\left( 1 + \frac{\Delta}{f_w(\bx)} \right)^{\gamma} -1 \; \right| \; \xi_t
=\bx \right],\]
and as long as $\Delta/f_w(\bx)$ is not too large we can use the fact that 
for $\gamma \in \R$ and small $x$ 
\[
(1+x)^{\gamma} = 1 + \gamma x + \frac{1}{2} \gamma (\gamma-1) x^2 + O(x^3).\]
Under the conditions of the lemma, for $\bx \in \W(\alpha)$ with $r = \|\bx\|$ large enough, a.s.,
\begin{align*}
 \left| \frac{\Delta}{f_w (\bx)} \right| \leq \frac{ C r^{w-1}}{f_w (\bx)} = O(r^{-1}),\end{align*}
using (\ref{fbnd}) and (\ref{lbd}).
 Hence for $\gamma \in \R$ and all $\|\bx\|$ large enough
\begin{align}
\label{ccc}
 & \Exp [ f_w(\xi_{t+1})^{\gamma} - f_w(\xi_t)^{\gamma} \mid \xi_t = \bx] 
= \gamma f_w (\bx) ^{\gamma-1} \Exp [ \Delta  \mid  \xi_t = \bx ] \nonumber\\
&~~~ {}+ \frac{1}{2} \gamma (\gamma-1) f_w (\bx)^{\gamma-2} \Exp [ \Delta^2 \mid \xi_t =\bx ] 
 + O \left( f_w (\bx)^{\gamma-3} \Exp [ \Delta^3 \mid \xi_t = \bx] \right).\end{align}
 Then from (\ref{ccc}), Lemma \ref{fmoms}, and (\ref{lbd})
  we obtain the desired result. 
\qed\\
 
 We will also need the following straightforward result.
 
 \begin{lemma}
 \label{777}
 Let $h: \R^2 \to \R$ and $R \subset \R^2$ be such that
$h(\bx) \leq 0$ for all $\bx \in \R^2 \setminus R$. Set
$\hat h (\bx) := h (\bx) \1_{\{ \bx \in R\}}$.
Then for  all $\bx \in R$ and all $t \in \Z^+$,
\[ 
 \hat h (\xi_{t+1} ) - \hat h (\xi_t) 
 \geq   h (\xi_{t+1} ) -   h (\xi_t), ~\textrm{on}~\{ \xi_t =\bx\}
.\]
\end{lemma}
\proof 
For $\bx \in R$ we have on $\{ \xi_t = \bx\}$ that 
$\hat h (\xi_{t+1} ) - \hat h (\xi_t)  
  = 
    h (\xi_{t+1} ) -   h (\xi_t)  
  -  h(\xi_{t+1}) \1_{\{\xi_{t+1} \notin R\}}$,
  which yields the result given that $h(\bx) \leq 0$ for $\bx \notin R$.   
  \qed

\section{Critical case: proofs of Theorems \ref{thm1}, \ref{thm2}, and \ref{thm9}}
\label{prf1}

\subsection{Overview and statement of upper bound}

In this section we prove our main results on moments of $\tau_\alpha$ in the critical
case, Theorems \ref{thm2} and \ref{thm9}, as well as giving a self-contained proof of Theorem \ref{thm1}
including the case $\alpha = \pi$ not directly covered by the results of \cite{mmw1}.
There are two largely separate components to the proofs of these three theorems.
The existence of moments part of Theorem \ref{thm2}, as well as Theorem \ref{thm1},
will follow from Lemma \ref{taillem} stated below, which gives
an upper bound on the tails of $\tau_\alpha$. We are not able to use the general results
such as Lemma \ref{aimthm} to prove Lemma \ref{taillem}; instead our proof is in some sense more
elementary, although we do use semimartingale tools at several points. On the other hand, for the
non-existence part of Theorem \ref{thm2}, as well as Theorem \ref{thm9}, we are able to appeal
 to the general result Lemma \ref{aimthm2} after finding and analysing a suitable Lyapunov function. Thus in the
second (non-existence) part of the proof the intuition is encapsulated in the Lyapunov function
and 
there is not a natural central lemma to stand alongside Lemma \ref{taillem} in that case.

Here is our central result for the `existence' part of the proofs. 

 \begin{lemma}
 \label{taillem}
 Suppose that (A1), (A2), and (\ref{drift1}) hold.
 Let $\alpha \in (0,\pi/2)$ and $\bx \in \W(\alpha)$. 
 There exist $\gamma \in (0,\infty)$, not depending on $\bx$,
 and $C \in (0,\infty)$, which does depend on $\bx$,
 such that for all $t>0$,
\begin{equation}
 \label{tails}
 \Pr [ \tau_\alpha > t \mid \xi_0 = \bx ] \leq C t^{- \gamma } .\end{equation}
 \end{lemma}

Now we describe the outline of the remainder of this section.
First, in Section \ref{sec:thm1prf}, we show how Lemma \ref{taillem} gives an almost immediate proof of Theorem \ref{thm1},
including the $\alpha = \pi$ case not covered by \cite{mmw1}. Crucial to the proof of Lemma \ref{taillem}
will be a {\em decomposition}
of the random walk $\Xi$ based
on the regularity condition (A1). In \cite[Section 4.2]{mmw1} we used a related decomposition that was, however, different, and in fact 
 more complicated
than the one used below, since  \cite{mmw1} considers general dimensions.
The version of the decomposition  in the present paper is described in detail in Section \ref{sec:deco}.
Section \ref{sec:rect} is devoted to a key step in the proof of Lemma \ref{taillem}, which is a result on the exit
from rectangles (Lemma \ref{rectangle}  below) that says, loosely speaking,
 that
if the walk starts somewhere near the centre of a rectangle,
there is strictly positive probability (uniformly
in the size of the rectangle) that the walk
will first exit the rectangle via the top/bottom. Here the fact that
$\|\mu(\bx)\| = O(\| \bx \|^{-1})$ is crucial. This result clarifies the key difference
between the one-dimensional and multi-dimensional settings: see the remark after Lemma \ref{rectangle}.
In Section \ref{sec:cones} we give the proof of Lemma \ref{taillem}. Then we turn to the `non-existence' parts
of the proof; our main tool is a 
 Lyapunov function introduced in Section \ref{sec:nonex}.
Finally we complete the proofs of Theorems \ref{thm2} and \ref{thm9} in Section \ref{sec:prfs}.

\subsection{Proof of Theorem \ref{thm1}}
\label{sec:thm1prf}

We establish Theorem \ref{thm1} 
 by studying the behaviour of the walk on a
set of seven overlapping quarter-planes that together
span $\W(\pi)$ (the plane minus a thickened half-line).
For this reason, 
we need to consider wedges like $\W(\alpha)$ 
with several different principal axes. This requires
some more notation.
Define lattice vectors $\bq_i$,
$i \in \{1,\ldots,7\}$ by
\[ \bq_5 = -\bq_1 = \be_1 + \be_2, \quad \bq_6 = -\bq_2  = \be_2, \quad \bq_3 =
-\bq_7 = \be_1 - \be_2, \quad \bq_4 = \be_1. \]
We also need
notation for 
perpendiculars to the $\bq_i$, specifically 
\[ \bq_i^{\perp} = \bq_{i+2},\ i\in \{ 1,2, \ldots, 5\}, \quad \bq_6^{\perp} =
-\bq_4, \quad \bq_7^{\perp} = \bq_1.  \]
For the corresponding unit vectors, write
$\hat \bq_i := \| \bq_i \|^{-1} \bq_i$ and $\hat \bq^\perp_i := \| \bq^\perp_i \|^{-1}
\bq_i^\perp$;
note that $\| \bq_i \| = \| \bq_i^\perp \|$, which is $1$ for even $i$ and $\sqrt{2}$
for odd $i$.
 
For $\beta \in (0,\pi/2)$ and
 $i \in \{ 1, \ldots, 7\} $
let $W_i(\beta)$ denote the wedge 
with apex $\0$,  internal angle $2\beta$, and
principal direction $\bq_i$; that is
\begin{equation}
\label{qdef}
W_i(\beta) := \{ \bx \in \R^2 : \bx \cdot \bq_i >0 , | \bx \cdot \bq_i^\perp | <  (\tan \beta) | \bx \cdot \bq_i |  \}.
\end{equation}
With our existing notation, this means that $W_4(\alpha)$ is $\W(\alpha)$;
the other $W_i(\alpha)$ are rotations of $\W(\alpha)$ through angles $k\pi/4$,
$k \in \{ \pm 1, \pm 2, \pm 3\}$.
In the proof of Theorem \ref{thm1} below
we will need  the quarter-planes $W_i(\pi/4)$;
when it comes to the
proof of Theorem \ref{thm2}
we need $W_i (\alpha)$ for $\alpha \in (0,\pi/2)$.
Thus we work in this generality for now.
For $\beta \in (0,\pi/2)$, let
 \begin{equation}
 \label{taubeta}
  \tau_i(\beta) := \min \{ t \in \Z^+ : \xi_t \notin W_i(\beta) \}. \end{equation}

\noindent
{\bf Proof of Theorem \ref{thm1}.}
It suffices to show that the result holds for $\alpha = \pi$, i.e., the walk a.s.\ 
eventually hits the thickened
half-line $\H_b$.
For notational ease let $Q_i := W_i (\pi/4)$,
$\tau_i := \tau_i (\pi/4)$ for $i \in \{1,\ldots,7\}$. Also write $Q_8 := \H_b$ and
 $B := B_A (\0)$ for some $A \in (0,\infty)$.
 
 Suppose that $\xi_0 \in Q_i$.
It follows immediately from Lemma \ref{taillem}
that $\Pr [ \tau_i < \infty ] =1$, and
 so $\Xi$
almost surely exits the initial quadrant $Q_i$. 
By the bounded jumps
assumption (A2), the definition (\ref{qdef}), and an appropriate choice of $A$,
we see that at time $\tau_i$, $\xi_{\tau_i}$ is
either: (i) in $Q_8$; (ii) in $B$; or (iii) within
distance $b$ of the principal axis of either $Q_{i+1}$
or $Q_{i-1}$, working mod $8$ for the indices of the $Q_j$s.

In case (ii) or (iii), $\Xi$ exits $B$ or the quadrant whose principal axis it is close to in finite time.
In the first case, having left $B$, $\Xi$ is in some quadrant, which it must exit in finite time, again ending
up in $B$ or close to the principal axis of some other quadrant.
This process repeats, showing that $\Xi$ must, infinitely often, be close to the principal axis of one
of the quadrants $Q_i$. Moreover, at such times, the proof of Lemma
\ref{lm4.1} shows that the events that the walk next visits  $Q_{i \pm 1}$ each have uniformly positive probability.
It follows that $\Xi$ visits each $Q_i$ eventually, a.s., and in particular hits the thickened half-line.
\qed 

\subsection{Decomposition}
\label{sec:deco}

For each $i$,
using the regularity condition (A1) we decompose $\Xi$
into a symmetric walk in the $\bq_i^{\perp}$ direction and
a residual walk. For $\bx, \by \in \Z^2$, $n \in \N$ and $t \in \Z^+$
let  
$p(\bx, \by;  n) : = \Pr [ \xi_{t+n} = \by \mid \xi_t = \bx ]$.
It follows from (A1) by considering finite combinations of jumps that
for each $i \in \{ 1, 2, \ldots, 7\}$ there exist constants $\gamma_i
\in (0,1/2)$, $n_i, j_i \in \N$ such that 
\begin{equation}
\label{0508a}
 \min_{\bx \in Q_i} \{ p(\bx,\bx + j_i \bq_i^{\perp}; n_i), 
 p(\bx, \bx -
j_i \bq_i^{\perp}; n_i) \} \geq \gamma_i  . \end{equation}

Now we fix $i$ and consider the `$n_i$-skeleton' random walk, i.e.~the
embedded
process $(\xi_{t n_i})_{t \in \Z^+}$.
For notational convenience, for $t \in \Z^+$ write
$\xi^*_t := \xi_{t n_i}$.
Then $\Xi^* = (\xi^*_t)_{t \in \Z^+}$
is a Markov chain on $\SS$
with transition probabilities
$\Pr [ \xi^*_{t+1} = \by \mid \xi^*_t = \bx ] = p( \bx, \by ; n_i )$, 
and $\xi^*_0 = \xi_0$.
 The walk $\Xi^*$
 inherits regularity from $\Xi$ as described in the following lemma.
 
 \begin{lemma} Suppose that (A1) and (A2) hold. Then
  \begin{align}
 \label{0508b}
  \Pr [ \| \xi^*_{t+1} - \xi^*_t \| \leq b n_i ] & = 1 ; ~\textrm{and}\\
  \label{minvar}
   \Exp [ | ( \xi^*_{t+1} - \xi^*_t ) \cdot \hat \bq_i^\perp |^2 \mid \xi^*_t = \bx ]
 &  \geq 2 j_i^2 \| \bq_i^\perp \|^2 \gamma_i >0 , \end{align}
 for all $\bx \in \SS$.
 Moreover, if (\ref{drift1}) holds, then, for $\bx \in \SS$, as $\| \bx \| \to\infty$,
\begin{equation}
\label{0530a}
 \| \Exp [ \xi^*_{t+1} - \xi^*_t \mid \xi^*_t = \bx ] \| = O(\| \bx \|^{-1} ) .\end{equation}
\end{lemma}
\proof 
The bound (\ref{0508b}) is immediate from (A2), while (\ref{minvar})
follows from (\ref{0508a}). Moreover, it follows from (A2) that,
\begin{equation}
\label{567}
 \max_{t n_i \leq s \leq (t+1)n_i}  \| \xi_s - \xi^*_t \| \leq n_i b, ~ {\rm a.s.}, \end{equation}
which with
 (\ref{drift1})
implies (\ref{0530a}). \qed\\

  By
(\ref{0508a}), 
 there exist sequences of random
variables $(V_t)_{t\in \N}$ and $(\zeta_t)_{t \in \N}$
such that:
\begin{itemize}
\item[(i)] 
  the $(V_t)_{t \in \N}$ are i.i.d.~with
 $V_t \in \{-1,0, 1\}$,
 $\Pr[V_t = 0] = 1-2\gamma_i$, and
 $\Pr [ V_t = -1 ] = \Pr [ V_t = +1 ] = \gamma_i$;
 \item[(ii)]
  $\zeta_{t+1} \in \Z^2$ with
 $\Pr [ \zeta_{t+1} =  \0 \mid V_t \neq 0 ] =1$; and
 \item[(iii)]
 we can decompose
 the jumps of $\Xi^*$ via, for $t \in \Z^+$,
\begin{equation}
\label{decomp}
 \xi^*_{t+1} - \xi^*_t = 
 \xi_{(t+1)n_i} - \xi_{tn_i} = \sum_{s=0}^{n_i - 1}
\theta(\xi_{tn_i +s}) =   V_{t+1} j_i \bq_i^{\perp} +  \zeta_{t+1} .  
\end{equation}
\end{itemize}
Note that given point (ii), (\ref{decomp}) is equivalent to, for $t \in \Z^+$,
 \[  \xi^*_{t+1} - \xi^*_t = 
   V_{t+1} j_i \bq_i^{\perp} \1_{\{ V_{t+1} \neq 0\}} +  \zeta_{t+1}  \1_{\{ V_{t+1} = 0\}} .   \]
 Thus we decompose the jump of $\Xi^*$ at time $t$ into a symmetric component
 in the perpendicular direction ($V_{t+1} j_i \bq_i^\perp$), 
 and a residual component ($\zeta_{t+1}$),  such that
 at any time $t$  only one of the two components is present in a particular
 realization. By (\ref{decomp}),
\begin{equation}
\label{xidec}
 \xi^*_t = \xi_0 + \sum_{s=1}^t ( V_s j_i \bq_i^\perp + \zeta_s ) .\end{equation}
This decomposition is valid throughout $\Z^2$, but for our purposes
we will apply the decomposition involving $\bq_i^\perp$ in the wedge
$W_i(\beta)$ for appropriate $\beta \in (0,\pi/2)$.

\subsection{Exit from rectangles}
\label{sec:rect}

We will use the decomposition of Section \ref{sec:deco} to establish (in Lemma
\ref{rectangle} below) how the walk exits from
sufficiently large rectangles aligned in the $\bq_i, \bq_i^\perp$ directions.
First we need two lemmas that deal in turn with the two parts of the
decomposition.

The rough outline of the proof of Lemma \ref{rectangle} below
is as follows.
In time $\lfloor \eps N^2 \rfloor$, we show that the process
driven by $V_1, V_2, \ldots$ will with positive
probability attain distance
sufficient to take it well beyond the top/bottom
of the rectangle; this is Lemma \ref{lm4.1} below.
On the other hand, we show that in time
$\lfloor \eps N^2 \rfloor$, for small enough $\eps>0$,
the residual process does not stray very far from its initial
point with good probability, regardless
of the realization of $V_1, V_2,\ldots$; this
is Lemma \ref{lm4.2} below. Together,
these two results will enable us to conclude
that with good probability
the walk will leave a rectangle via the top/bottom. 
First we need some more notation. 
Set $Y_0 :=    \xi_0 \cdot \hat \bq_i^\perp $ and, for $t \in \N$,
\begin{equation}
\label{ytdef}
 Y_t := Y_0 + j_i \| \bq^\perp_i \| \sum_{s=1}^{t}   V_s . \end{equation}
Then $Y_t$
 is 
 the displacement of the symmetric part
of the decomposition for $\xi^*_t$ in the
 $\bq_i^{\perp}$ direction. The process $( Y_t
)_{t \in \Z^+}$ is a symmetric, homogeneous
 random walk on $\| \bq^\perp_i \| \Z$ with
$\Pr [ Y_t = Y_{t-1} ] =\Pr [ V_t =0] = 1 - 2\gamma_i < 1$ and 
jumps of size $\| \bq^\perp_i \|j_i$. For $h \in (0,\infty)$, let
\begin{equation}
\label{0508e}
  \tau^\perp_h := \min \left\{ t\in \Z^+ : |Y_t| \geq \lceil 3h N\rceil \|  \bq_i^{\perp} \| 
\right\}.  \end{equation} 
 
\begin{lemma}
\label{lm4.1}
 Suppose that  (A1) holds. 
Let $h \in (0,\infty)$.
For any $\eps>0$, there exist $\delta>0$ and $N_1 \in \N$ such that
for any $N \geq N_1$ and any $y  \in \Z$ with $|y| \leq 2 h N$,
\[ \Pr [ \tau^\perp_h \leq \lfloor \eps N^2\rfloor \mid Y_0 =  \|  \bq_i^{\perp} \| y ] \geq  \delta . \] 
\end{lemma}
\proof Fix $h \in (0,\infty)$.
Suppose $Y_0 =  \|  \bq_i^{\perp} \| y$, $|y| \leq 2h N$.
 If $y \neq 0$ then couple a copy of the walk $Y_t$ started
from $ \|  \bq_i^{\perp} \| y$ with another $\tilde{Y}_t$ started from 
$0$ which has jumps in the opposite direction to $Y_t$ until $| Y_t - \tilde
Y_t| \leq  \|  \bq_i^{\perp} \| j_i$ for the first time, 
from which time on $Y_t, \tilde Y_t$
jump in the same direction.
Then when $| \tilde Y_t | \geq K$ we have $|Y_t| \geq K -  \|  \bq_i^{\perp} \| j_i$,
and with probability $\gamma_i$ the next jump will take $| Y_{t+1} | \geq K$.
It follows that for any $\eps>0$,
\begin{align*} \Pr [ \tau^\perp_h \leq \lfloor \eps N^2 \rfloor \mid Y_0 =  \|  \bq_i^{\perp} \| y ]
& \geq \gamma_i \Pr [ \tau^\perp_h \leq
\lfloor \eps N^2 \rfloor -1 \mid Y_0 = 0 ] \\
& \geq \gamma_i \Pr [ \tau^\perp_h \leq \lfloor \eps' N^2\rfloor \mid Y_0 = 0 ] ,\end{align*}
for any $\eps' \in (0,\eps)$ and all $N$ large enough. 
 Hence it suffices to take $y = 0$.

The process $(Y_t)_{t \in \Z^+}$ is a symmetric random
walk on $ \|  \bq_i^{\perp} \| \Z$ with
independent, bounded jumps and $\Exp[|Y_{t+1} - Y_t|^2] = 2  \gamma_i \|  \bq_i^{\perp} \|^2 j_i^2 > 0$.
Standard 
central limit theorem estimates imply
that for any $\eps>0$ there exists $\delta>0$
such that for all $N$ sufficiently large
\[ 
 \Pr \left[  Y_{\lfloor \eps N^2 \rfloor} > 
 \lceil 3h N \rceil \| \bq_i^{\perp}\|  \mid Y_0 = 0 \right] \geq \delta ,
\textrm { and }
 \Pr \left[  Y_{\lfloor \eps N^2 \rfloor} < 
 -\lceil 3h N \rceil  \|  \bq_i^{\perp}\|  \mid Y_0 = 0 \right] \geq \delta .\]
Each of the (disjoint) events in the last
display implies that $\tau^\perp_h \leq \lfloor \eps N^2\rfloor$.
 \qed
\\

Let $Z_0 : = (\xi_0 \cdot \hat \bq_i ) \hat \bq_i$ and for
$t \in \N$ let 
\begin{equation}
\label{ztdef}
 Z_t := Z_0 + \sum_{s=1}^{t}  \zeta_s .\end{equation} 
Thus $(Z_t)_{t \in \Z^+}$
is the residual part of the process $(\xi^*_t)_{t \in \Z^+}$
after the symmetric perpendicular process $(Y_t)_{t \in \Z^+}$
has been extracted. Indeed, with $Y_t, Z_t$ as defined
at (\ref{ytdef}), (\ref{ztdef}) we have  $\xi^*_0 = \xi_0 = Y_0  \hat \bq_i^\perp
+ Z_0$, and also from (\ref{xidec}) that for $t \in \N$,
\begin{equation}
\label{decomp2}
 \xi^*_t = Y_t \hat \bq_i^\perp + Z_t .\end{equation}

We next show that with good probability the residual process
$(Z_t)_{t \in \Z^+}$ does not exit from a suitable ball around its
initial point by time $\lfloor \eps N^2\rfloor$. By construction 
  $(Z_t)_{t \in \Z^+}$ depends  on $  ( V_t) _{t \in \N}$
since  the distribution of $\zeta_{t+1}$ depends  on  $\xi_{tn_i}$. For $t \in \N$,
 let $\Omega_V(t) := 
  \{-1,0, 1\}^{t}$ and let $\omega_V
\in \Omega_V(t)$ denote a generic
realization of   $( V_1, \ldots, V_t)$.
 
\begin{lemma}
\label{lm4.2}
 Suppose that  (A1), (A2), and (\ref{drift1}) hold. 
Let $r \in (0,1/2]$.
There exist $N_2 \in \N$ and $\eps>0$ such that
for all $N \geq N_2$, all $z  \in \Z$ with $| z | \leq b$,
and all $\omega_V \in \Omega_V ( \lfloor \eps N^2 \rfloor)$,
\[ \Pr \left[ \max_{0 \leq t \leq \lfloor \eps N^2 \rfloor} 
\| Z_t - Z_0 \| \leq
r N \mid (V_1,\ldots, V_{\lfloor \eps N^2 \rfloor} ) = \omega_V, Z_0 = (N+z ) \bq_i    \right] \geq \frac{1}{2} .
\] 
\end{lemma}
\proof
Although the decomposition used in the present paper is different, the proof of this result
is similar to (in fact, due to the stronger regularity assumptions, simpler than) the proof
of the corresponding Lemma 4.5 in \cite{mmw1}, so we omit it. \qed 
\\

We now define notation
for our rectangles. Fix $h \in (0,\infty)$,
which will determine the aspect ratio
of the rectangles. For $N \in \N$, let 
\begin{equation}
\label{0508c}
 S(N) := \{ \bx \in \Z^2 : 0 < \bx \cdot \hat \bq_i < 2 N\| \bq_i \|,
   |  \bx \cdot \hat \bq_i^\perp | < 2h N\| \bq_i\| \} ,
 \end{equation}
and also define regions adjacent to $S(N)$ via
\begin{align}
\label{0508d}
 U_1(N) & :=  \{ \bx \in \Z^2 : \bx \cdot \hat \bq_i \geq  2 N \| \bq_i \|\},\nonumber\\
  U_2(N) & := \{ \bx \in
\Z^2 : 0 < \bx \cdot \hat \bq_i < 2N \| \bq_i\|, | \bx \cdot \hat \bq_i^\perp| \geq 2h N\| \bq_i\| \}  .    
\end{align} 
Lemmas \ref{lm4.1} and \ref{lm4.2}
 combine to enable us show
 that $\Xi$ exits $S(N)$ via $U_2(N)$ with good probability
when started from somewhere near the   bisector of $S(N)$ in the $\bq_i^\perp$
direction. 

\begin{lemma}
 \label{rectangle}
 Suppose that  (A1), (A2), and (\ref{drift1}) hold. 
 Let $h \in (0,\infty)$.
 Then there exist $\delta > 0$, $N_0 \in \N$
  such that
 for any $N \geq N_0$, any $y, z \in \Z$ with $| y| \leq 2h N$ and $|z| \leq b$,  
 \[ \Pr [ \Xi \textrm{ hits } U_2(N) \textrm{ before } U_1(N) \mid \xi_0 =
 (N+z) \bq_i + y \bq_i^\perp ] \geq \delta . \]
\end{lemma}

\rem
This result highlights the   difference between the 
exit-from-cones problem   and the analogous
problem of exit from a half-line in one-dimension, where drift $O(x^{-1})$
does {\em not}  imply finiteness of the exit time. The one-dimensional
analogue of Lemma \ref{rectangle} is false:
 classical  
gambler's ruin estimates imply that for a  random walk
on $\Z^+$ with mean-drift $O(x^{-1})$ at $x$, the probabilities
of hitting $0$, $2M$ first, starting from $M$,
 are not necessarily bounded uniformly away from $0$. \\
 
 \noindent
 {\bf Proof of Lemma \ref{rectangle}.}
Fix $h \in (0,\infty)$.
Suppose that $\xi_0 =
 (N+z ) \bq_i + y \bq_i^\perp$.
 Let $\eps>0$ be as in the $r = (1\wedge h)/2$
 case of Lemma \ref{lm4.2}.
 Suppose that $N \geq \max \{ N_1,N_2\}$ with $N_1, N_2$
as in Lemmas \ref{lm4.1}, \ref{lm4.2} respectively.
Define the events
\[ G := \left\{ \max_{0 \leq t \leq \lfloor \eps N^2 \rfloor} \| Z_t - Z_0 \| \leq (1 \wedge h) N/2 \right\}, ~~~
H := \left\{\tau^\perp_h \leq \lfloor \eps N^2 \rfloor \right\}. \]
 By (\ref{decomp2}) we have that $| \xi^*_t \cdot \hat \bq_i^\perp | = | Y_t + Z_t \cdot  \hat \bq_i^\perp |
= |Y_t + (Z_t -Z_0) \cdot \hat \bq_i^\perp |$,
since $Z_0 \cdot \hat \bq_i^\perp =0$. It follows by the
triangle inequality
that on
  $G \cap H$, 
\[ | \xi^*_t \cdot \hat \bq_i^\perp | \geq | Y_t |  
- \| Z_t -Z_0 \| \geq \lceil 3h N \rceil \| \bq^\perp_i \| - (1 \wedge h) (N/2) \geq 2h N \| \bq_i^\perp \| ,\]
for some $t \leq \lfloor \eps N^2 \rfloor$, which in particular
implies that $| \xi_t \cdot \hat \bq_i^\perp | \geq 2h N \| \bq^\perp_i \|$ 
for some $t \leq n_i \lfloor\eps N^2\rfloor$.
On the other hand, also on $G \cap H$ it follows
from (\ref{decomp2})  that
\begin{align*}
 \max_{0 \leq t \leq n_i \lfloor \eps N^2 \rfloor} | \xi_t \cdot \hat \bq_i |
\leq \max_{0 \leq t \leq \lfloor \eps N^2 \rfloor} | \xi^*_t \cdot \hat \bq_i | + n_i b
= \max_{0 \leq t \leq \lfloor \eps N^2 \rfloor} | Z_t \cdot \hat \bq_i | + n_i b \\
\leq | Z_0 \cdot \hat \bq_i | + \max_{0 \leq t \leq \lfloor \eps N^2 \rfloor} \| Z_t - Z_0 \|
+ n_i b < 2N \| \bq_i^\perp \| ,\end{align*}
for all $N$ sufficiently large, since $Z_0 \cdot \hat \bq_i =
\xi_0 \cdot \hat \bq_i = (N+z ) \| \bq_i\|$.
Hence (with $\xi_0$ as given)
\[ E: = \left\{ \Xi \textrm{ hits } U_2(N) \textrm{ before } U_1(N) \right\} \supseteq G \cap H .\]
$H$ is determined by the realization
$\omega_V \in \Omega_V (\lfloor \eps N^2 \rfloor)$, and  so (with $\xi_0$ as given)
\[ \Pr [ E   ] \geq
\Pr[ G \cap H   ]
= \sum_{\omega_V \in \Omega_V (\lfloor \eps N^2 \rfloor ) : H ~{\rm occurs}}
  \Pr[ G
\mid \omega_V  ] \Pr[  \omega_V    ] .\]
Applying Lemma \ref{lm4.2} with $r = (1 \wedge h)/2$ to $\Pr [ G \mid \omega_V ]$ we then obtain
 \begin{align*}
 \Pr[ E \mid \xi_0 = (N+z ) \bq_i + y \bq_i^\perp ] \geq  
\frac{1}{2} \sum_{\omega_V \in \Omega_V (\lfloor \eps N^2 \rfloor ) : H ~{\rm occurs} }
  \Pr[    \omega_V    ] 
  = \frac{1}{2} \Pr [ H   ] \geq \frac{\delta}{2} > 0 ,
  \end{align*}
 applying Lemma \ref{lm4.1}.
     \qed
 
 \subsection{Exit from cones: proof of Lemma \ref{taillem}}
 \label{sec:cones}
 
 We can now prove our key upper tail bound.
  Recall the definition of $\tau_i(\beta)$ from (\ref{taubeta}).\\
 
 \noindent
 {\bf Proof of Lemma \ref{taillem}.}
 Take $\xi_0  \in W_i (\beta)$, $\beta \in (0,\pi/2)$.
Let $h = \tan \beta \in (0,\infty)$ and
 \[
k_0 := \min \{ k \in \N : 2^k \| \bq_i \| \geq
 \xi_0 \cdot \hat \bq_i , 2^k \geq N_0, 2^k \geq b\} ,\]
 where $N_0$ is as in Lemma \ref{rectangle}
 and $b$ is as in (A2).
 Consider the sequence of rectangles $S (2^k)$
where $k \in \Z^+$, as defined at (\ref{0508c}), with $h=\tan \beta$.
 Set $\sigma_0 :=0$ and, for $k \in \N$,
\[ \sigma_k := \min \{ t \in \Z^+ : \xi_t \cdot \hat \bq_i \geq 2^k \| \bq_i\| \} .\]

Suppose that $\Xi$ has not left $W_i (\beta)$ by the time $\sigma_k$
for some $k \geq k_0$, i.e.,
$\tau_i(\beta) > \sigma_k$.
Then, using (A2),
$2^k \| \bq_i \| \leq \xi_{\sigma_k} \cdot \hat \bq_i \leq 2^k \| \bq_i \| +b$
and on $\{ \tau_i(\beta) > \sigma_k \}$, from (\ref{qdef}),
\[  
|\xi_{\sigma_k} \cdot\hat \bq^\perp_i | < h | \xi_{\sigma_k} \cdot \hat \bq_i |
 \leq 2^k h \| \bq_i \| + hb \leq  2 \cdot 2^k h  \| \bq_i \| ,\]
  for all $k \geq k_0$. Hence we can apply
Lemma \ref{rectangle} to the walk started
at $\xi_{\sigma_k}$, with $N=2^k \geq N_0$
for $k \geq k_0$. Then, with $U_1(N), U_2(N)$
as defined in (\ref{0508d}), we obtain, for all $k \geq k_0$,
\begin{equation}
\label{0508f}
 \Pr[ (\xi_t)_{t \geq \sigma_k} \textrm{ hits } U_2 (2^k) \textrm{ before } U_1 (2^k)
\mid \tau_i(\beta) > \sigma_k ] \geq \delta >0.
 \end{equation}
 But by definition of $U_1(N), U_2(N)$, and (\ref{qdef}),
 we have that
if $\Xi$ hits $U_2 (N)$ before $U_1 (N)$, then
$\Xi$ leaves the wedge $W_i(\beta)$. 
Moreover, if $\Xi$ has not hit $U_1(2^k)$
by time $\tau_i(\beta)$, then 
$\max_{0 \leq t \leq \tau_i(\beta)} \xi_t \cdot \hat \bq_i < 2^{k+1} \| \bq_i \|$,
so that $\tau_i(\beta) < \sigma_{k+1}$. 
Hence the inequality (\ref{0508f})
can   be expressed as $\Pr [ \tau_i(\beta) \leq \sigma_{k+1} \mid \tau_i(\beta) > \sigma_k ] \geq \delta >0$,
for all $k \geq k_0$. Hence, for all $k > k_0$,
\begin{align}
\label{0509a}
\Pr[ \tau_i(\beta) > \sigma_k ] = 
\prod_{j=k_0+1}^k \Pr [ \tau_i(\beta) > \sigma_j \mid \tau_i(\beta) >\sigma_{j-1} ] \cdot \Pr[ \tau_i(\beta) > \sigma_{k_0} ]
\leq C (1-\delta)^k ,
\end{align}
for some $C=C(k_0,\delta) \in (0,\infty)$ that does not depend on $k$.

We next  estimate the tails of the times $\sigma_k$.
It is most convenient to work once again
via the embedded walk $\Xi^*$. Set $\sigma^*_k := \min \{ t \in \Z^+ : \xi^*_t \cdot  \hat \bq_i \geq 2^k \| \bq_i \| \}$.
For $t \in \Z^+$, for the remainder of this proof write
 $X_t := \xi^*_t \cdot \hat \bq_i$. Let $A,C>0$
and set
$W_t := ((C+X_{t \wedge \tau_i (\beta)})^A)$.
We show that for $A,C$ sufficiently large,
the process  $(W_t)_{t \in \Z^+}$  is a strict submartingale so that we can apply
a result from \cite{mvw} to obtain a bound for $\Exp [ \tau_i(\beta) \wedge \sigma^*_k]$.

Note that Taylor's theorem implies that for any $x \geq 0$ and any $y \in \R$
 with $|y|$ bounded,
\begin{align*}
 (C+x+y)^A - (C+x)^A 
= A (C+x)^{A-1} \left[ y + \frac{(A-1) y^2}{2 (C+x)} + O  ((C+x)^{-2} ) \right] .\end{align*}
Set  $\theta_t^* = \xi^*_{t+1} - \xi^*_t$. 
Let $\F_t = \sigma (\xi_0, \ldots \xi_t)$. 
By (\ref{0508b})
we may apply the last displayed
equation 
with $x = \xi_t^* \cdot \hat \bq_i$ and $y = \theta^*_t \cdot \hat \bq_i$ and take expectations
to obtain 
\begin{align*}
& \Exp [ W_{t+1} - W_t  \mid \F_{n_i t} ] 
 \\
{} = {} & A (C+ X_t )^{A-1} 
 \left[    \Exp [ \theta_t^* \cdot \hat \bq_i \mid \F_{n_i t} ] 
+ \frac{ (A-1)}{2} \frac{ \Exp [ ( \theta_t^* \cdot \hat \bq_i )^2 \mid \F_{n_i t} ]}{C+ X_t}
+ O ( (C+X_t )^{-2} ) \right], \end{align*}
on $\{ t < \tau_i (\beta) \}$. Also, on $\{ t < \tau_i (\beta) \}$ we have that $X_t \geq 0$ and
$X_t \leq \| \xi^*_t \| \leq O (X_t)$. So
 using (\ref{minvar}) and (\ref{0530a}) we have that the last displayed expression is
 bounded below by
 \[  A (C+ X_t )^{A-1} 
 \left[    - C_1 (1+ X_t)^{-1}   
+ \frac{ (A-1)}{C_2} (C+ X_t)^{-1}
+ O ( (C+X_t )^{-2} ) \right] ,\]
for some constants $C_1, C_2 \in (0,\infty)$. Hence 
we can choose $A, C$ sufficiently large so that
$\Exp [ W_{t+1} - W_t  \mid \F_{n_i t} ] \geq \eps > 0$, 
on $\{ t < \tau_i (\beta) \}$. Moreover, from (\ref{0508b}) we have
$| X_{t+1} - X_t | \leq \| \xi^*_{t+1} - \xi^*_t \|   \leq   n_i b$.
Hence we can apply a straightforward modification of \cite[Lemma 3.2]{mvw} 
to obtain, for all $k \geq k_0$,
$\Exp [ \tau_i (\beta) \wedge \sigma^*_k ]  \leq
  \eps^{-1} (C+ 2^{k+1} +  n_i b)^A$. 
By definition of $\xi^*_t$,  $\sigma_k \leq n_i \sigma^*_k$ a.s.,
hence  there exists $C \in (0,\infty)$
such that $\Exp [ \tau_i (\beta) \wedge \sigma_k ] \leq  2^{kC}$,
for all $k \geq k_0$.
Markov's inequality with $M=C+1$ then implies that for $k \geq k_0$,
\[ \Pr [ \tau_i (\beta) \wedge \sigma_k > 2^{kM} ] \leq 2^{-kM} \Exp [ \tau_i (\beta) \wedge \sigma_k] 
\leq 2^{kC} \cdot 2^{-kM} = 2^{-k} .\]
Combining this with (\ref{0509a}) and the fact that for any $k$,
\begin{align*}
\Pr [ \tau_i(\beta) > 2^{kM} ] \leq \Pr [ \tau_i(\beta) > \sigma_k ]
+ \Pr [ \tau_i (\beta) \wedge \sigma_k > 2^{kM} ] , \end{align*}
we have that 
$\Pr [ \tau_i(\beta) > 2^{kM} ] \leq C  (1-\delta)^k + 2^{-k}$
for all $k \geq k_0$.
 It follows that
 there exist constants $M, \gamma' \in (0,\infty)$, not depending on $\bx$,
 and $C \in (0,\infty)$, which does depend on $\bx$,
 such that for all $k \geq k_0$,
\begin{equation}
 \label{tailss}
 \Pr [ \tau_i(\beta) > 2^{kM} \mid \xi_0 = \bx ] \leq C 2^{- \gamma' k} .\end{equation}
Clearly the result extends to all $k \in \Z^+$
for a suitable choice of $C$ in (\ref{tailss}),
depending on $k_0$ and so also on $\xi_0$.
 For any $t >0$, we have that
 $t \in [ 2^{kM}, 2^{(k+1)M} )$
 for some $k \in \Z^+$. Then given $\xi_0 = \bx$ we have
 from (\ref{tailss}) that
 \[ \Pr [ \tau_i (\beta) >t ]
 \leq \Pr [ \tau_i (\beta) > 2^{kM} ]
 \leq C 2^{-\gamma' k}
 \leq C ( 2^{-M} t)^{-\gamma'/M} = C'  t^{-\gamma} ,\]
 for some $C', \gamma \in (0,\infty)$,
 not depending on $t$, with, moreover, $\gamma$  
 not depending
 on $\bx$.
  \qed
  
  \subsection{Non-existence of moments via an almost-linear Lyapunov function}
  \label{sec:nonex}

 This section
is devoted to our technique for
establishing
the non-existence part of the proof
of Theorem \ref{thm2}, which will also
enable us to give
 a proof of Theorem \ref{thm9}.
In the wedge $\W(\alpha)$, $\alpha  \in (0,\pi/2)$,
 we work 
with the embedded walk $\xi^*_t = \xi_{t n_i}$, where
in this case we can take $n_i = n_0$ as in (A1).
 We first aim to
show that for any $\alpha \in (0,\pi/2)$ there
exists $p \in (0,\infty)$ such that $\Exp [ \tau_\alpha^p ] =\infty$.

The outline of our approach is  as follows.
We consider a one-dimensional process $(Y_t)_{t \in \Z^+}$ where
$Y_t = g(\xi^*_t)$ for 
a suitably chosen
$g$ and apply Lemma \ref{aimthm2}.
More specifically, we construct
an almost linear or $\eps$-linear
(in the sense of Malyshev \cite{malyshev}, see also \cite{kingman} and \cite[Chapter 3]{fmm})
function $g$ to enable us to apply the generalized
form \cite{aim} of ``Lamperti's
conditions'' \cite{lamp2}
in Lemma \ref{aimthm2}.
 The idea is to construct $g$ so that its
level curves are  horizontal translates of
$\partial\mathcal{W}(\alpha)$ but with the 
apex replaced by a circle arc.

\begin{figure}[!h]
\centering
\includegraphics[width=13cm]{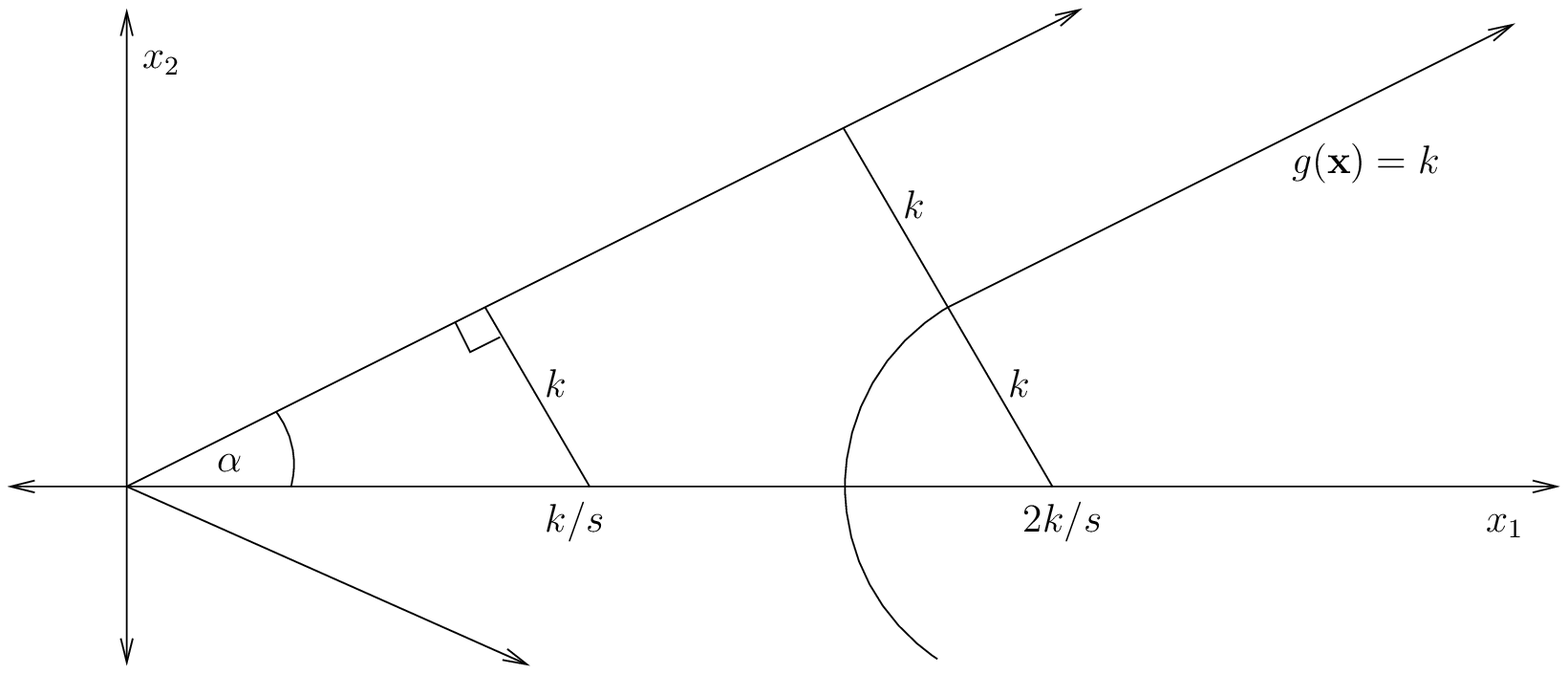}
\caption{Level curve of the function $g$.}
\label{fig2}
\end{figure}

Fix $\alpha \in (0,\pi/2)$.
During the remainder of this section,
set $s := \sin \alpha \in (0,1)$, $c := \cos \alpha \in (0,1)$. 
We now construct the function $g: \R^2 \to [0,\infty)$. 
Set $g(\bx) = 0$ for $\bx \in \R^2 \setminus \W(\alpha)$.
For $\bx =(x_1,x_2)\in \W(\alpha)$ such that
$|x_2| \geq  \frac{sc}{1+c^2} x_1$ let
$g(\bx) = s x_1 - c | x_2|$.
For $\bx \in \W(\alpha)$ with
$|x_2| \leq \frac{sc}{1+c^2} x_1$, set
$g(\bx) = k \in [0,\infty)$
on the minor arc of the circle
\begin{equation}
\label{circle}
((2k/s)-x_1)^2 + x_2^2 = k^2
\end{equation}
 between
$(k(1+c^2)/s,kc)$ and
$(k(1+c^2)/s,-kc)$.
Then $g$ is specified on $\W(\alpha)$ by its level
curves $g(\bx) = k$, $k \geq 0$,
each of which is $\partial \W(\alpha)$ 
translated so that the apex is at $(k/s,0)$
and the tip of the wedge
smoothed to a circular arc. See Figure \ref{fig2}.

We now state some properties of the function $g$.
Observe that for $\bx \in \R^2$,
\begin{equation}
\label{gbnd}
g(\bx) \leq \| \bx \|. \end{equation}
For $\bx \in \W(\alpha)$ with $|x_2| \geq \frac{sc}{1+c^2} x_1$,
 $\nabla g(\bx) = (s, \pm c)$ and $\|
\nabla g(\bx) \| =1$; 
for $|x_2| \leq \frac{sc}{1+c^2} x_1$,
\begin{equation}
\label{0531b}
 \nabla g(\bx) = \frac{1}{D(\bx)} ( - ((2g(\bx)/s)-x_1), x_2 ) = \frac{1}{D(\bx)} 
 \left( - \sqrt{ g(\bx)^2 -x_2^2 } , x_2 \right)
,\end{equation}
from (\ref{circle}), where
$D(\bx) := g(\bx)+ (2/s) (x_1 - (2g(\bx)/s))$. 
When $|x_2| \leq \frac{sc}{1+c^2} x_1$, so that the level curve of $g$ is 
a circular arc, we have
\begin{align}
\label{arc}
g(\bx) ((2/s)-1) & \leq x_1 \leq  g(\bx) ((2/s)-s), ~\textrm{and} \\
\label{dbounds}
-g(\bx) ((2/s)-1) & \leq D(\bx) \leq - g(\bx). \end{align}
It follows from (\ref{0531b})
that for $\bx \in \W(\alpha)$ with
$|x_2| \leq \frac{sc}{1+c^2} x_1$,
$\| \nabla g(\bx) \|  = | D(\bx) |^{-1} g (\bx)$.
Hence from (\ref{dbounds}), 
\begin{align}
\label{grad}
\inf_{\bx \in \W(\alpha)} \| \nabla g(\bx) \| \geq \frac{s}{2-s} \geq \frac{s}{2} ,
~\textrm{and}~
\sup_{\bx \in \R^2} 
\| \nabla g(\bx) \| \leq 1 .
\end{align}
 
To obtain our non-existence of moments result
for $\tau_\alpha$, we will apply Lemma \ref{aimthm2}
to  $Y_t = g(\xi^*_t)$. 
The next  lemma 
gives some further
 properties of  $g$ that we will need
 here and later in the proof of Theorem \ref{thm9}. 

\begin{lemma}
For $\bx \in \W(\alpha)$ with $|x_2| \leq \frac{sc}{1+c^2} x_1$
we have
\begin{equation}
\label{0601b}
 g(\bx) \geq \frac{s}{2} \| \bx \| .\end{equation}
Also, there exists $\eps>0$ such that for all $\bx \in \W(\alpha)$,
\begin{equation}
\label{0806c}
D_1 g (\bx) \geq \eps .\end{equation}
Finally, there exists $C \in (0,\infty)$ such that
for all $\bx \in \R^2$ and all $i,j \in \{1,2\}$,
\begin{equation}
\label{0601a}
 | D_{ij} g (\bx) | \leq C \| \bx \|^{-1} .\end{equation}
\end{lemma}
\proof
To obtain (\ref{0601b}), we observe that 
for $|x_2| \leq \frac{sc}{1+c^2} x_1$, from
 (\ref{arc}),
 \begin{align*} 
\| \bx \|^2 =
x_1^2 + x_2^2 \leq \left[ \left( \frac{sc}{1+c^2} \left( (2/s) -s \right)\right)^2
+ \left( (2/s) -s \right)^2 \right] g(\bx)^2 \\
= \left[ c^2 + ((2/s)-s)^2 \right] g(\bx)^2 = ((4/s^2)-3) g(\bx)^2 ,\end{align*}
and (\ref{0601b}) follows. 
Consider   (\ref{0806c}).
It suffices to suppose $|x_2| \leq \frac{sc}{1+c^2} x_1$.
By (\ref{0531b}),
\begin{equation}
\label{0601c}
 D_1 g(\bx) = \frac{2g(\bx) - sx_1}{((4/s)-s) g(\bx) - 2x_1}
= R \left[ 1 + \frac{Sx_1}{g(\bx)-Rx_1} \right] ,\end{equation}
where $R \in (0,2/3)$ and $S \in (0,1/6)$ are defined as
\begin{equation}
\label{RSdef}
 R = \frac{2}{(4/s)-s}, ~ \textrm{and}~ S = R - (s/2) = \frac{s}{2} \left( \frac{s^2/4}{1-(s^2/4)} \right).\end{equation}
Here, since $s \in (0,1)$, it is straightforward to show that in fact
\begin{equation}
\label{Rbound}
 \frac{s}{2} \leq R \leq \frac{2s}{3}, ~ \textrm{and} ~ \frac{s^3}{8} \leq S \leq \frac{s^3}{6}. \end{equation}
Moreover, we have from (\ref{RSdef}) and (\ref{arc}) that
\begin{equation}
\label{0601d}
 (s^2/4) g(\bx)  \leq g(\bx)-Rx_1  \leq (s/2) g(\bx) .\end{equation}
It follows from (\ref{0601c}) and (\ref{0601d})
that $D_1 g(\bx) \geq R$, and so with (\ref{Rbound}) 
 we get (\ref{0806c}).

Now consider (\ref{0601a}).
Note that $D_{ij} g(\bx) = 0$ unless
$\bx \in \W(\alpha)$ with $|x_2| \leq \frac{sc}{1+c^2} x_1$,
so it suffices to consider that case.
First   consider $D_{11} g(\bx)$.
Differentiating in (\ref{0601c}) yields
\begin{align}
\label{0601e}
D_{11} g(\bx) & = \frac{RS}{g(\bx) - R x_1}
- \frac{RSx_1}{(g(\bx)-Rx_1)^2} \left( D_1 g(\bx) - R \right)
\nonumber\\
& = \frac{RS}{(g(\bx)-Rx_1)^3} \left( (g(\bx)-Rx_1)^2 - RS x_1^2 \right),
\end{align}
using (\ref{0601c}) once more.
Then from (\ref{0601e}), using   (\ref{Rbound}), (\ref{0601d}), and (\ref{arc}), 
together
with (\ref{0601b}),
we obtain (\ref{0601a}) in the case $i=j=1$. 
The other cases of (\ref{0601a}) follow by analogous
but tedious calculations, 
which we omit.
 \qed

 \subsection{Proofs of Theorems \ref{thm2} and \ref{thm9}}
\label{sec:prfs}

The next result 
gives some basic properties of the process $(g(\xi^*_t))_{t \in \Z^+}$.

\begin{lemma}
\label{lem0530}
  Suppose that (A1), (A2) and (\ref{drift1}) hold.
  Then there exist $B, C \in (0,\infty)$ and $\eps >0$ for which, for any $\bx \in \W(\alpha)$,
  \begin{align}
  \label{0531a}
 \Pr [  | g(\xi^*_{t+1} ) - g(\xi^*_t ) | \leq B ] & =1 ; \\
    \label{0530b}
   |\Exp[g(\xi^*_{t+1}) - g(\xi^*_t) \mid \xi^*_t =
  \bx ]| & \leq C \| \bx \|^{-1}  ; \\
  \label{0531c}
  \Exp[ (g(\xi^*_{t+1}) - g(\xi^*_t))^2 \mid \xi^*_t =
  \bx ] & \geq \eps .\end{align}
\end{lemma}
\proof The mean value theorem for functions of two variables
 implies that 
\begin{equation}
\label{0531d}
g(\xi^*_{t+1}) - g (\xi^*_t) = ( \xi^*_{t+1} - \xi^*_t ) \cdot \nabla g ( \bz ) ,
\end{equation}
where $\bz = \xi^*_t + \eta ( \xi^*_{t+1}-\xi^*_t)$ for some $\eta \in [0,1]$.
So (\ref{0531d}) implies that $| g(\xi^*_{t+1}) - g(\xi^*_t) | \leq \| \xi^*_{t+1} - \xi^*_t \|$,
a.s., by (\ref{grad}), which with (\ref{0508b}) yields (\ref{0531a}).
Similarly, by (\ref{0531d}) and (\ref{grad}),
\[ | \Exp [ g(\xi^*_{t+1}) - g(\xi^*_t) \mid \xi^*_t = \bx ] |
\leq 2 \| \Exp [ \xi^*_{t+1} - \xi^*_t \mid \xi^*_t = \bx] \| ,\]
and then 
(\ref{0530a}) implies (\ref{0530b}). 
Finally, using (A1) we have from (\ref{0531d}) that
for $\bx \in \W(\alpha)$,
\begin{align*} \Exp[ (g(\xi^*_{t+1}) - g(\xi^*_t))^2 \mid \xi^*_t =
  \bx ] \geq \kappa
  [ k_1 \be_1 \cdot \nabla g(\bz_1) ]^2  
  ,\end{align*}
where $\bz_1 = \bx + \eta_1 k_1 \be_1$,
for some $\eta_1 \in [0,1]$;
so in particular $\bz_1 \in \W(\alpha)$.
Thus, by (\ref{0806c}),
$\Exp[ (g(\xi^*_{t+1}) - g(\xi^*_t))^2 \mid \xi^*_t =
  \bx ]
\geq \kappa k_1^2  \eps^2
> 0$, giving (\ref{0531c}).
\qed\\

Now we verify that $g(\xi^*_t)$
satisfies the conditions of Lemma \ref{aimthm2}.

\begin{lemma}
\label{lem0603}
Suppose that (A1), (A2) and (\ref{drift1}) hold.  
For $A>0$ large enough
 there exist $C, D \in (0,\infty)$, $r>1$, and $p_0>0$ such that for any $t \in \Z^+$,
on $\{ \upsilon_A > t\}$, (\ref{noncon1}), (\ref{noncon2}), and (\ref{noncon3}) hold
for $Y_t = g(\xi^*_t)$.
\end{lemma}
\proof
Let $Y_t = g(\xi^*_t)$, $t \in \Z^+$. Let $r > 0$.
We need to estimate $\Exp [ Y_{t+1}^{2r} - Y_t^{2r} \mid \xi^*_t = \bx]$.
By Taylor's theorem, 
for $y >0$ and $\delta$ with $|\delta | \leq B$, there
exists $\eta  \in [0,1]$
for which
\begin{align}
\label{y2r}
(y+\delta)^{2r} - y^{2r} 
& = 2r \delta y^{2r-1} + r (2r-1) \delta^2 (y + \eta \delta)^{2r-2}\nonumber \\
& = 2r \delta y^{2r-1} + r (2r-1) \delta^2 y^{2r-2} + o(y^{2r-2}) .\end{align}
We now establish (\ref{noncon2}). Let $r=1$ in (\ref{y2r}) to obtain
\[ \Exp [ Y_{t+1}^2 - Y_t^2 \mid \xi^*_t = \bx ]
\geq 2 g(\bx) \Exp [ Y_{t+1} - Y_t \mid \xi^*_t = \bx ]
 \geq -2C g(\bx) \| \bx \|^{-1} ,\]
 by (\ref{0530b}). Then (\ref{gbnd})
 completes the proof of (\ref{noncon2}).
 Now let $\F_t = \sigma ( \xi_0^*, \ldots , \xi_t^*)$.
Then, by the $r>1$ case of
(\ref{y2r}),
$\Exp [ Y_{t+1}^{2r} - Y_t ^{2r} \mid \F_t ]$ is bounded above by 
\[ 2r g(\bx)^{2r-1} \Exp [ Y_{t+1} - Y_t \mid \F_t ]
+ 2 r^2 \Exp [ (Y_{t+1} - Y_t)^2 \mid \F_t ] (Y_t + B)^{2r-2} .\]
On $\{ \upsilon_A > t \}$, $g (\xi^*_t ) > A$ so $\xi^*_t \in \W(\alpha)$.
So by (\ref{0531a}) and (\ref{0530b}), on $\{ \upsilon_A > t\}$,
\[ \Exp [ Y_{t+1}^{2r} - Y_t ^{2r} \mid \F_t ]
 \leq   2C r Y_t^{2r-1}  \| \xi^*_t \|^{-1} + 2 r^2 B^2 (Y_t + B)^{2r-2} = O ( Y_t^{2r-2} ),\]
 by (\ref{gbnd}). Thus (\ref{noncon3})
is satisfied for $r>1$. Similarly, 
taking $r = p_0$ in (\ref{y2r}) and using (\ref{0530b}) again, but this time using
the lower bound in
(\ref{0531c}), valid on $\{ \upsilon_A > t\}$,
\begin{align*}
\Exp [ Y_{t+1}^{2p_0} - Y_t^{2p_0} \mid \F_t ]
& \geq  - 2 p_0 C Y_t^{2r-1} \| \xi^*_t \|^{-1} + p_0 ( 2p_0-1) \eps Y_t^{2r -2}  + o ( Y_t^{2r-2} ) \\
& \geq Y_t^{2r-2} p_0 \left( -2C + (2p_0 - 1) \eps + o(1) \right), \end{align*}
by (\ref{gbnd}), and the last expression
 is non-negative 
on $\{ \upsilon_A >t\}$,
taking $A$ and $p_0$ sufficiently large. \qed\\

 We need one more result that says, under our regularity conditions,
 an asymptotically zero drift ensures that the walk cannot be forced to jump
 straight out of the wedge with probability 1, provided it starts far enough away from $\0$.
 
 \begin{lemma}
 \label{lem88}
 Suppose that (A1) and (A2) hold and that $\| \mu (\bx) \| \to 0$ as $\| \bx \| \to \infty$.
 Let $\alpha \in (0,\pi]$.
 There exist $\eps, A,C \in (0,\infty)$ such that
  for any $\bx \in \W(\alpha)$ with $\| \bx \| \geq A$,
 \[ \Pr [ \Xi ~\textrm{hits}~ B_C ( ( \eps \|\bx \|, 0 ) ) ~ \textrm{before} ~ \R^2 \setminus \W(\alpha) \mid \xi_0 = \bx ] > 0 .\]
 \end{lemma}
 \proof 
 Let $d(\bx)$ denote the distance of $\bx$ from the boundary of the wedge $\W(\alpha)$.
 Suppose that $\bx \in \W(\alpha)$ and, without loss of generality, $x_2 >0$.
 First let $\alpha < \pi/2$.
 Given $d(\bx) > b n_0$, condition (A1) implies that with probability at least $\kappa$
 the walk starting at $\bx \in \W(\alpha)$ will end up at $\bx - k \be_2$ in $n_0$ steps, while
 during this time (A2) implies the walk cannot have left the wedge. Repeating this argument a finite
 number of times (depending on $\bx$)
 until the desired ball is reached leads to the desired conclusion for all such $\bx$.
 A similar argument works when $\alpha \geq \pi/2$ and $d(\bx) > bn_0$, starting with  steps of $k \be_1$.
 
 Thus it remains to deal with the case where the walk starts at $\bx$ with $d (\bx) \leq b n_0$ but $\| \bx \|$ large.  
 Recall (see \cite[p.\ 4]{mvw}) that (A1) implies that $\Pr [ \xi_{t+1} \neq \bx \mid \xi_t = \bx]$
 is uniformly positive. 
 We may suppose that
 $\bx$ is such that $\Pr [ (\xi_{t+1} -\xi_t) \cdot \be_\perp (\alpha) \neq 0 \mid \xi_t = \bx] >0$,
 since if this is not the case then (A1) entails that there is positive probability of the walk reaching such an $\bx$
 in a finite number of jumps parallel to the boundary of the wedge (and hence, by (A2), remaining inside $\W(\alpha)$ provided the walk
 started far enough from $\0$).
 So we may take $\bx$ such that there is positive probability of the next jump
 having a component perpendicular to the wedge boundary. In fact, for $\|\bx\|$ large enough, we have
 $\Pr [ (\xi_{t+1} -\xi_t) \cdot \be_\perp (\alpha) < 0 \mid \xi_t = \bx] >0$, so that there is
 positive probability of the walk jumping `farther into' the wedge.
 To see this, note that since $\Xi$ lives on (a subset of) $\Z^2$
 and, by (A2), has uniformly bounded jumps
 there are only finitely many possible values for $\xi_{t+1} - \xi_t$, and so any  non-zero
 component in the $\be_\perp (\alpha)$ direction must in fact be greater in absolute
 value than some $\delta >0$ not depending on $\bx$. Then
 \[ \mu (\bx ) \cdot \be_\perp (\alpha)
  \geq \delta \Pr [ (\xi_{t+1} - \xi_t ) \cdot \be_\perp (\alpha) \geq \delta
   \mid \xi_t = \bx ] - b \Pr [ (\xi_{t+1} - \xi_t ) \cdot \be_\perp (\alpha) < \delta
   \mid \xi_t = \bx ]  .\]
 Take $\| \bx \|$ large enough so that
 $\| \mu (\bx) \| \leq \eps$. Writing $p = \Pr [ (\xi_{t+1} - \xi_t ) \cdot \be_\perp \geq \delta \mid \xi_t = \bx ]$
 we have
 $\eps \geq \delta p - b(1-p)$, implying that $p < 1$
 for $\eps$ small enough. For $\|\bx\|$ large enough, a finite number of such jumps occur with positive
 probability and take $\Xi$ to distance at least $b n_0$ from the boundary of the
 wedge, so we can then appeal to the first part of the proof.
  \qed\\

\noindent
{\bf Proof of Theorem \ref{thm2}.}
 Let $\alpha \in (0,\pi/2)$.
 To prove Theorem \ref{thm2},
 it suffices to show that
 $\Exp[ \tau_\alpha^p ] =\infty$
 and $\Exp [ \tau_\alpha^q] < \infty$
 for some $p,q$ with
 $0<q<p<\infty$.
First, we apply 
 Lemma \ref{taillem}
 in the case $i=4$, $\beta = \alpha$,
 so that $W_4 (\beta) = \W(\alpha)$
 and $\tau_i(\beta)= \tau_\alpha$.
   Then  
 from Lemma \ref{taillem}, for some $\gamma, C \in (0,\infty)$,
   where  $\gamma$ does not
 depend on $\bx$,
\[ \Exp [ \tau^q \mid \xi_0 = \bx ] \leq 1 + \int_1^\infty \Pr [ \tau > r^{1/q} ] \ud r
\leq 1+C \int_1^\infty r^{-\gamma/q} \ud r < \infty ,\]
provided $q \in (0,\gamma')$.

Finally, 
Lemma \ref{lem0603}
implies that we can apply Lemma
\ref{aimthm2}
with $Y_t = g (\xi_t^*)$, so
that for some $A,p \in (0,\infty)$ we have
$\Exp [ \upsilon_A  ^p \mid \xi_0 = \bx] =\infty$
for all $\bx \in \W(\alpha)$ with  $g(\bx)$  sufficiently large.
But by definition of $g$ and $\Xi^*$, and (A2),
$\tau_\alpha \geq n_0 (\upsilon_A -1)$,  a.s.,
for $A>b$. Hence $\Exp [ \tau_\alpha ^p \mid \xi_0 = \bx] =\infty$
for all $\bx \in \W(\alpha)$ with $g (\bx)$
  sufficiently large. By Lemma \ref{lem88}, the conclusion extends to all $\bx \in \W(\alpha)$
  with $\| \bx \|$ large enough.
 \qed\\

\rem The difficulty with extending the overlapping
quadrant argument of Section \ref{sec:thm1prf}
to show existence of moments for $\alpha \geq \pi/2$
is that the constant $C$ in Lemma \ref{taillem}
depends upon $\bx$, and so some control is required
over the location of $\Xi$ on its exit
from each quadrant $Q_i$. Such an argument
should be possible using similar techniques to those employed here;
for reasons of space we do not pursue this here.\\

To prove Theorem \ref{thm9}, we use a similar argument to the non-existence
part of Theorem \ref{thm2}.
In particular,
we refine the lower bound in (\ref{0530b})
that depends explicitly upon the constant
$c$ in (\ref{qq2}). 
For this (in Lemma \ref{lem0601} below),
we replace the first-order
Taylor expansion used in the proof
of Lemma \ref{lem0530} with a second-order
expansion. 
   
 \begin{lemma}
 \label{lem0601}
 Suppose that (A1) and (A2) hold, and that $\alpha \in (0,\pi/2)$.
 Suppose that (\ref{qq2})
 holds for some $c>0$. Then there
 exist $\eps, C \in (0,\infty)$, not depending on $d$, such that
 for all $\bx \in \W(\alpha)$ with $\| \bx \|$ sufficiently large
 \[ \Exp [ g( \xi_{t+1}^*) - g(\xi_t^* ) \mid \xi^*_t = \bx]
 \geq  \| \bx \|^{-1} ( \eps c - C ) .\]
\end{lemma}
\proof Write $\xi^*_{t+1} - \xi^*_t
= ( \theta^*_1 (\xi^*_t), \theta^*_2 (\xi^*_t) )$
in Cartesian components. 
Given $\xi^*_t = \bx$,  (\ref{567}) holds with $n_i = n_0$, so 
(\ref{qq2}) implies that $\Exp [ \theta^*_1 (\bx) ] \geq (n_0 c + o(1) ) \| \bx \|^{-1}$
and $\Exp [ \theta^*_2 (\bx) ] = o ( \| \bx \|^{-1} )$.
Conditional on $\xi^*_t = \bx$,
 Taylor's theorem gives
\[ g(\xi_{t+1}^*) - g(\xi_t^*)
= \sum_i \theta_i^* (\bx) D_i g(\bx)
+ \frac{1}{2} \sum_{i,j} \theta_i^* (\bx) \theta_j^* (\bx) D_{ij} g (\bz) \]
for some $\bz \in \R^2$.
Then taking expectations and using
 (\ref{grad}), (\ref{0601a}), and (A2),
 we have
 \[ \Exp  [ g( \xi_{t+1}^*) - g(\xi_t^* ) \mid \xi^*_t = \bx]
 \geq \frac{n_0 c +o(1)}{\| \bx \|} D_1 g(\bx) - \frac{C }{\| \bx\|}.\]
 Then (\ref{0806c}) 
 completes the proof. \qed\\
 
 \noindent
 {\bf Proof of Theorem \ref{thm9}.}
 Again we  apply Lemma \ref{aimthm2} to $Y_t = g (\xi_t)$,
 analogously to the proof
  of the non-existence part of Theorem
 \ref{thm2}. Repeating the argument for Lemma \ref{lem0603}, (\ref{noncon2}) and (\ref{noncon3}) hold as before,
 but now using 
 Lemma \ref{lem0601} we have that (\ref{noncon1}) holds for $p_0 =1/2$, taking $c$ sufficiently large.
\qed

\section{Subcritical case: proof of Theorem \ref{thm5}}
\label{prf2}

\subsection{Overview}

 In this section we prove Theorem \ref{thm5}.
 The proofs of parts (i) and (ii) of Theorem \ref{thm5}
use the   Lyapunov
  functions $f_w$ defined at (\ref{fdef}) but are otherwise independent.
  
The functions $f_w$ are well-suited to the subcritical case,
allowing us to obtain the explicit exponents in Theorem \ref{thm5}.
 The technique in Section \ref{sec:cones},
 used to prove the existence of moments in
 Theorem \ref{thm2}, does not give
 sharp exponents, since $\gamma$ in Lemma \ref{taillem}
 depends on the $\delta$ in Lemma \ref{lm4.1}, which depends
 upon the $\eps$ in Lemma \ref{lm4.2}, and these
 results assume very general conditions on $\Xi$. In addition,
 the method in Section \ref{sec:cones}  works only for $\alpha < \pi/2$.
Similarly, the method used in Section
\ref{sec:nonex} to prove the non-existence of moments
in Theorem \ref{thm2} is not sharp enough to produce
the correct exponents
that we require
for Theorem \ref{thm5}, and again needs $\alpha < \pi/2$.

The outline of this section is as follows. In Section \ref{exmom} we
give the technical preliminaries for the proof of Theorem \ref{thm5}(i)
which we present in Section \ref{sec:exprf}. The more difficult
problem of non-existence of moments  needs considerably more work.
Preliminary calculations are in Section \ref{nonexmom}. We are not able to use the general result
Lemma \ref{aimthm2} in this case, so we use a more elementary approach based on giving a lower bound
for the probability that the walk takes a certain time to leave a wedge. This key estimate is given in Section
\ref{sec:key}. Finally, the proof of Theorem \ref{thm5}(ii) is completed in Section \ref{sec:prfnonex}.

\subsection{Existence of moments}
\label{exmom}
  
  For a given $\alpha \in (0,\pi]$,
 we fix $w \in (0,\pi/(2\alpha))$. Then $\W(\alpha)$
 lies inside the larger wedge $\W(\pi/(2w))$.
   Define the modified random walk $\tilde \Xi = (\tilde \xi_t)_{ t \in \Z^+}$ 
  by $ \tilde \xi_t := \xi_t \1_{ \{ t \leq \tau_\alpha \}}$, so that $\tilde \Xi$
is identical to $\Xi$
 on $\W(\alpha)$ but from $\bx \notin \W(\alpha)$ jumps directly to $\0$ and remains
 there; then $\tilde \xi_t = \0$ for $t \geq \tau_\alpha+ 1$.
For $t \in \Z^+$, set $X_t := f_w (\tilde \xi_t ) ^{1/w}$.
For $B \in (0,\infty)$, define
\[  
\tilde \tau_{\alpha,B}   :=
\min \{ t \in \Z^+ : X_t \leq B   \} .\]
The next result will be the basis for our results in this section.

\begin{lemma}
\label{exlem}
Suppose that (A1) and (A2) hold.
Fix $\alpha  \in (0,\pi]$ and $w \in (0,\pi/(2\alpha))$.
Suppose that there exist $p_0 >0$,
 $A_0,C \in (0,\infty)$
such that for all $\bx \in \W(\alpha)$ with $\| \bx \| \geq A_0$,
\begin{align}
\label{eee}
\Exp [ f_w(\xi_{t+1})^{2p_0/w}  - f_w (\xi_t) ^{2p_0/w} \mid \xi_t = \bx ] 
\leq - C f_w (\bx) ^{(2p_0-2)/w} .\end{align}
Then for  any $p \in [0,p_0)$ and any
  $\bx \in \W(\alpha)$,
$\Exp [ \tau_{\alpha}  ^p \mid \xi_0 = \bx] < \infty$.
\end{lemma}
\proof 
 $\tilde \xi_{\tau_\alpha +1} = \0$ so  
 $X_{\tau_\alpha +1} =0$; hence, for any $B>0$,
 $\tilde \tau_{\alpha, B} \leq \tau_\alpha +1$ a.s.. Hence
\begin{align}
\label{0621a}
 \{ \tilde \tau_{\alpha, B} > t \}
\subseteq 
\{ \tau_\alpha > t, \xi_t \in \W(\alpha) \}
\cup \{ \tau_\alpha = t, \xi_t \in \W(\pi/(2w)) \setminus \W(\alpha) \} ,
\end{align}
for all $B$ sufficiently large,
using (A2) and the fact that by definition $\{ \tilde
\tau_{\alpha,B} > t\} \subseteq \{ \| \xi_t \| > B \}$.
We consider the two events in the disjoint union in (\ref{0621a}) in turn.
Let $\F_t := \sigma ( \xi_0, \xi_1 ,\ldots , \xi_t )$.
 On $\{ \tau_\alpha>t\}$ we have $\xi_t  = \tilde \xi_t$
 and $\xi_{t+1} = \tilde \xi_{t+1}$. So by
(\ref{eee})
  there exists $C' \in (0,\infty)$
such that,
on $\{ \tau_\alpha > t\}$,
\begin{align}
\label{ngl}
\Exp [ X_{t+1}^{2p_0} - X_t^{2p_0} \mid \F_t ] \leq - C' X_t^{2p_0-2}.\end{align}
On $\{ \tau_\alpha = t\}$, 
$ \Exp [ X^{2p_0}_{t+1} - X_t^{2p_0} \mid \F_t ] =  -  X_t^{2p_0}$, so that
on $\{ \tilde \tau_{\alpha,B}   > t \} \cap \{ \tau_\alpha = t\}$,
$ \Exp [ X^{2p_0}_{t+1} - X_t^{2p_0} \mid \F_t ] 
\leq - B^2 X_t^{2p_0-2}$,
since $\tilde \tau_{\alpha,B} > t$ implies that $X_t^2 \geq B^2$.

Thus,  for some $C'\in (0,\infty)$, (\ref{ngl})
holds  on $\{\tilde \tau_{\alpha,B}   > t\}$
for any $B \geq B_0$, say.
We apply Lemma \ref{aimthm} with $Y_t = X_t$
to obtain,  for any $p \in [0,p_0)$, $B \geq B_0$, and $\bx \in \W (\alpha)$,
\begin{align}
\label{mom2}
\Exp [ \tilde \tau_{\alpha,B}   ^ p \mid \xi_0 = \bx ] < \infty.\end{align}
It remains to
deduce the corresponding result for $\tau_\alpha$.

On $\{ \tau_\alpha \geq t \}$, $\tilde \xi_t = \xi_t$ and so by (\ref{lbd}),
on $\{ \tau_\alpha \geq t \}$, 
\begin{align}
\label{hhh} \|  \xi_t \| \geq X_t \geq \eps_{\alpha,w}^{1/w} \|  \xi_t\| 
.\end{align}
Recall that $\tilde \tau_{\alpha,B} \leq \tau_\alpha +1$.
On $\{\tilde \tau_{\alpha, B} \leq \tau_\alpha \}$,
$\| \xi_{\tilde \tau_{\alpha,B}  } \| \leq \eps_{\alpha,w}^{-1/w} X_{\tilde
\tau_{\alpha,B} } \leq \eps_{\alpha,w}^{-1/w} B$ by   
 (\ref{hhh}).
On the other hand, on
$\{\tilde \tau_{\alpha,B} = \tau_\alpha +1\}$,
clearly $\tau_\alpha \leq \tilde \tau_{\alpha, B}$.
 Recalling the definition of $\tau_{\alpha,A}$ from (\ref{taua}),
 it follows that for all $A \geq B \eps_{\alpha,w}^{-1/w}$, a.s.,
$\tau_{\alpha,A}    \leq  \tilde \tau_{\alpha,B}$.
Then with (\ref{mom2}) we obtain that
for all $\bx \in \W(\alpha)$  and all $A$ sufficiently large
$\Exp [   \tau_{\alpha,A}^ p \mid \xi_0 = \bx ] < \infty$.

Condition (A1) then extends the result
 to $\tau_\alpha$ by standard `irreducibility' arguments. Indeed,
 (A1) implies that for random variables $K_0, K_1, K_2, \ldots$
 with $\Pr [ K_i \geq t ] \leq \re^{-ct}$, for some $c >0$,
 $\tau_\alpha \leq \sum_{i=1}^{K_0} (\tau_i + K_i)$, where $\tau_1, \tau_2, \ldots$ are
 copies of $\tau_{\alpha, A}$; here $K_0$ represents the number of visits to $B_A(\0)$ before
 leaving the wedge, and $K_1, K_2, \ldots$ are the durations of the successive visits. By (A2),
 on each exit from $B_A(\0)$ into $\W(\alpha)$, $\Xi$ is restricted to a finite number of states, and
 so $\Exp [ \tau_i^p ]$ is uniformly bounded. Hence, for any $p < p_0$,
 \begin{align*} \Pr [ \tau_\alpha  \geq t ]  & \leq \Pr [ K_0 \geq C \log t ] + \sum_{i=1}^{C \log t} \Pr [ K_i \geq C \log t ]
 + \sum_{i=1}^{C \log t} \Pr [ \tau_i \geq t/(2C \log t) ] \\
 & = O ( t^{-p} (\log t)^{p+1} ) ,\end{align*}
 for $C < \infty$ large enough, by Boole's and Markov's inequalities. Thus $\Exp [\tau_\alpha^q ] < \infty$ for any $q < p < p_0$.
 This completes the proof.
 \qed

\subsection{Proof of Theorem \ref{thm5}(i)}
\label{sec:exprf}

The next result, with Lemma \ref{exlem},
will enable us to deduce Theorem \ref{thm5}(i).
 
 \begin{lemma}
 \label{dd2}
 Suppose that (A2) holds.
 Let $\alpha \in (0,\pi]$.
 Suppose that for some $\sigma^2 \in (0,\infty)$,
 for $\bx \in \W(\alpha)$
  as $\|\bx\| \to \infty$,
 \begin{equation}
 \label{c1}
  \| \mu (\bx) \|  = o(\| \bx\|^{-1});
 ~~~M_{12} (\bx) =o(1);
 ~~~M_{11} (\bx) = \sigma^2 +o(1); ~~~ M_{22} (\bx) = \sigma^2 + o(1).\end{equation}
 Then for any $w \in (0,\pi/(2\alpha))$ and any $\gamma \in (0,1)$,
 there exist constants $A,C \in (0,\infty)$ for which
    for all $\bx \in \W_A(\alpha)$,
 \begin{align}
 \label{bbbb}
  \Exp [ f_w (\xi_{t+1})^{\gamma} - f_w (\xi_{t})^\gamma 
 \mid \xi_t = \bx ] \leq - C f_w (\bx)^{\gamma-(2/w)}.\end{align}
\end{lemma}
\proof
Let $w \in (0, \pi/(2\alpha))$.
For $\bx \in \W(\alpha)$, we have that
(\ref{lbd}) holds. Then
by Lemma \ref{fgmoms}
with (\ref{c1})
 we have  that
 for $\gamma \in \R$,
  \begin{align}
  \label{ppp}   \Exp [ f_w(\xi_{t+1})^\gamma - f_w(\xi_t)^\gamma \mid \xi_t = \bx] 
 = \frac{1}{2} \gamma (\gamma-1) w^2 \sigma^2 f_w(\bx)^{\gamma-2} r^{2w-2} (1+o(1)), 
 \end{align}
 for all $\bx \in \W(\alpha)$, as $\| \bx \| \to \infty$. 
 It follows from (\ref{ppp}) and (\ref{lbd}) that 
for $\gamma \in (0,1)$ and some $C \in (0,\infty)$,
 for $\bx \in \W (\alpha)$ with $\| \bx \|$ sufficiently large (\ref{bbbb}) holds.
\qed\\
  
\noindent
{\bf Proof of Theorem \ref{thm5}(i).}
Let $w \in (0,\pi/(2\alpha))$.
For $\gamma \in (0,1)$,  take $p_0 = \gamma w/2$.
Then 
Lemma \ref{dd2}
says  that for $\bx \in \W(\alpha)$ 
with $\| \bx\|$ sufficiently large
(\ref{eee}) holds for $\gamma\in (0,1)$ and $w \in (0,\pi/(2\alpha))$.
Then Lemma \ref{exlem} implies that
for any $\bx \in \W(\alpha)$,
 $\Exp [ \tau_\alpha   ^p \mid \xi_0 = \bx ] < \infty$
for all $p \in [0,p_0)$. Since both $\gamma <1$ and $w < \pi/(2\alpha)$ may be
taken
 arbitrarily close to their upper bounds,
we may choose any $p$ less than $\pi/(4\alpha)$. \qed

 \subsection{Non-existence of moments}
 \label{nonexmom}

Let $\alpha \in (0,\pi]$. Throughout this section
we will take $w = \pi/(2\alpha)$. We will
again be interested in   $f_w (\xi_t)^\gamma$, $\gamma \in \R$,
this time in the wedge $\W(\alpha)$. 
Due to difficulties with estimating
the behaviour of $f_w(\xi_t)^\gamma$ near the
boundary of the wedge $\W(\alpha)$ (cf Lemma
\ref{fgmoms}), we cannot apply
the non-existence theorems from
\cite{aim} (such as Lemma \ref{aimthm2} above).
Thus we need a different approach.

A key step in this section is a good-probability lower-bound
on the time taken to leave a wedge; this is Lemma
\ref{lem2} below. A similar approach is used in
\cite{bfmp}, where Lemma 6.2 deals with a 
 special case
of a random walk in a quarter-plane. In any case, to show non-existence
of moments something like Lemma \ref{lem2} is required; analogous
lemmas are needed for the general results of \cite{lamp2,aim}.

We use the Lyapunov  function
$\hat f_w$ where 
$\hat f_w (\bx) := f_w (\bx) \1_{\{ \bx \in \W (\alpha) \}}$ for $\bx \in \R^2$.
The first task of this section is to estimate the mean increment of 
$\hat f_w (\xi_t)^\gamma$
  for  $\gamma>1$.
We recall that for $w = \pi/(2\alpha)$,
Lemma
\ref{fgmoms} applies for $f_w$ only in a  wedge smaller than $\W(\alpha)$. The next result
will allow us to
overcome this obstacle. For $K>0$ we use the notation
\begin{equation}
\label{wkdef}
 \W^K (\alpha) := \left\{ \bx \in \W(\alpha) : f_w (\bx) \geq K^{-1} \| \bx\|^{w-1} \right\}.
 \end{equation}

\begin{lemma}
\label{edge}
Let $\alpha \in (0,\pi]$
and $w=\pi/(2\alpha)$. 
Suppose that  (A2) holds, and that
for some $v \in (0,\infty)$,
 for $\bx \in \W(\alpha)$, as $\|\bx\| \to \infty$,
 \begin{equation}
 \label{sub2a} \| \mu(\bx) \| =o(1); ~~~
   M_{12} (\bx) =o(1);
 ~~~M_{11} (\bx) \geq v +o(1); ~~~ M_{22} (\bx) \geq v + o(1).\end{equation}
 Then there
  exist $A,K \in (0, \infty)$ such that
$ \Exp [ \hat f_w (\xi_{t+1}) - \hat f_w (\xi_t ) \mid \xi_t = \bx ] \geq 0$ 
for all $\bx \in \W_A (\alpha) \setminus \W^K (\alpha)$.
\end{lemma}
\proof
For   $K>0$, take $\bx \in \W(\alpha) \setminus \W^K (\alpha)$.
By (\ref{wkdef}),
$f_w (\bx) \leq K^{-1} r^{w-1}$ and hence
$\cos (w\varphi) \leq K^{-1} r^{-1}$.
Thus $\bx$ is  close to the boundary $\partial \W(\alpha)$. In order
to estimate the expected change in $f_w$ on a jump of $\Xi$ started
from $\bx$, we introduce the notation
$U(\bx) := \{ \by \in \W(\alpha) : f_w(\by) \geq f_w (\bx) \}$. 
   We use the shorthand $\hat \Delta = \hat f_w (\xi_{t+1}) - \hat f_w (\xi_t )$.

Since
$\hat f_w (\xi_{t+1}) \geq 0$, we have
$ \Exp [ \hat \Delta \1_{\{ \xi_{t+1} \notin U(\bx)\}} \mid \xi_t = \bx ]
\geq - f_w(\bx) \geq -K^{-1} r^{w-1}$, 
so  \begin{align}
  \label{xx1}
 \Exp [ \hat \Delta \mid \xi_t = \bx ] \geq
   \Exp [ \hat \Delta \1_{\{ \xi_{t+1} \in U(\bx)\}} \mid \xi_t = \bx ] 
  - K^{-1} r^{w-1}. \end{align}
For a random variable $X$ with  $\Pr[ |X| < m] =1$, 
$\Pr [ m|X| > X^2] =1$ and so $\Exp |X| \geq m^{-1} \Exp [X^2]$.
Moreover, 
 $\Exp [X \1_{\{ X \geq 0\}}] = (\Exp[X]+\Exp |X|)/2$. So we
conclude that 
\begin{equation}
\label{ew1}
\Exp [X \1_{\{ X \geq 0\}}] \geq (\Exp[X]+ m^{-1}\Exp[ X^2] )/ 2.
\end{equation}
Now write $\Delta = f_w (\xi_{t+1}) - f_w (\xi_t)$.
Then $\{ \Delta \geq 0, \xi_t = \bx\} = \{ \xi_{t+1} \in U (\bx), \xi_t = \bx \}$.
Hence
applying the elementary inequality (\ref{ew1}) with $X = \Delta$ 
and using the bound (\ref{fbnd}) gives, for some $C \in (0,\infty)$ and
all $\bx \in \W(\alpha)$,  
\begin{align*} & \Exp [   \Delta \1_{\{ \xi_{t+1} \in U(\bx)\}} \mid \xi_t = \bx ]  \\
& {} \geq  \frac{1}{2} \Exp [ f_w (\xi_{t+1}) -   f_w (\xi_t ) \mid \xi_t = \bx ]
+ C(1+ \| \bx \|)^{1-w} 
\Exp [ (f_w (\xi_{t+1}) -   f_w (\xi_t ) )^2 \mid \xi_t = \bx ].\end{align*}
By (\ref{sub2a}), we obtain from (\ref{fmom1})
and (\ref{fmom2}) 
 that
there exists $C >0$, not depending on $K$,
 such that,  for all $\bx \in \W(\alpha)$ with $\| \bx\|$ large enough,
\begin{equation}
\label{ii2}
 \Exp [  \Delta \1_{\{ \xi_{t+1} \in U(\bx)\}} \mid \xi_t = \bx ] 
\geq  C \| \bx \|^{w-1}. \end{equation}
It follows from Lemma \ref{777}
that we can replace $\Delta$ by $\hat \Delta$ in (\ref{ii2}).
Then the claimed result follows from (\ref{xx1})
 with (\ref{ii2}), by taking
$K$ large enough.
\qed\\

Here then is our result on the mean increment of  $\hat f_w (\xi_t)^\gamma$
  for $\gamma >1$.

 \begin{lemma}
 \label{sub1}
 Suppose that (A2) holds.
 Suppose that for
  some $\sigma^2 \in (0,\infty)$,
  as $\|\bx\| \to \infty$, (\ref{sub2}) holds.
 Let $\alpha \in (0,\pi]$.
 Then for $w = \pi/(2\alpha)$ and any $\gamma >1$,
 there exists $A \in (0,\infty)$ for which,
 for all $\bx \in \W_A(\alpha)$,
 \begin{align}
 \label{sub3}
  \Exp [ \hat f_w (\xi_{t+1})^{\gamma} - \hat f_w (\xi_{t})^\gamma 
 \mid \xi_t = \bx ] \geq 0 .\end{align}
\end{lemma}
\proof 
It suffices to take $\gamma \in (1,2]$.
Under the conditions of the lemma, Lemma \ref{edge} implies that for some $K$
the desired result holds
for $\bx \in \W_A (\alpha) \setminus \W^K (\alpha)$. So it
remains to consider $\bx \in \W_A(\alpha) \cap \W^K (\alpha)$.
Writing 
$\hat \Delta = \hat f_w (\xi_{t+1}) - \hat f_w (\xi_t )$,
we have
that for $\xi_t = \bx$,
\begin{equation}
\label{xx4}
 \hat f_w (\xi_{t+1})^\gamma
- \hat f_w (\xi_t)^\gamma =
(\hat f_w ( \bx ) + \hat \Delta)^\gamma - \hat f_w (\bx)^\gamma
= f_w (\bx)^\gamma   \left[ \left(1 + \frac{\hat \Delta}{\hat f_w (\bx)} \right)^\gamma -1 \right] .\end{equation}
To obtain a lower bound, we make use of the fact that for any $\gamma \in (1,2]$
and $L \in (0,\infty)$,
\begin{equation}
\label{xx3} (1+x)^\gamma \geq 1 + \gamma x + \frac{1}{2} (1+L)^{\gamma-2} \gamma (\gamma-1) x^2  \end{equation}
for $x \in [-1,L]$. To apply
(\ref{xx3}) in (\ref{xx4}) with $x = \hat \Delta / \hat f_w(\bx)$ we need
 $-\hat f_w (\bx) \leq \hat \Delta \leq L \hat f_w(\bx)$.
The first inequality here is automatically
satisfied since $\hat f_w (\xi_{t+1}) \geq 0$ a.s.. 
For the second inequality, we have  for $\bx \in \W^K (\alpha)$ from 
(\ref{fbnd}) and (\ref{wkdef}) that on $\{ \xi_t =  \bx\}$,
\[ \| \hat \Delta \|   \leq C \| \bx \|^{w-1} \leq C K f_w (\bx) = CK \hat f_w (\bx) .\]
So taking $L = C K$  we can indeed apply (\ref{xx3}) in (\ref{xx4}) to  obtain, for some $A, C \in (0,\infty)$,
for any $\bx \in \W_A(\alpha) \cap \W^K (\alpha)$,
conditional on $\xi_t = \bx$,
\[  \hat f_w (\xi_{t+1})^\gamma
- \hat f_w (\xi_t)^\gamma   \geq
\gamma f_w (\bx)^{\gamma -1} \hat \Delta + C f_w (\bx)^{\gamma-2} \hat \Delta^2.\]
The right-hand side of the last display   is increasing in $\hat \Delta$, and
so by Lemma \ref{777} we can replace $\hat \Delta$
by $\Delta$ and then take expectations
to obtain
\begin{align*} \Exp [ \hat f_w (\xi_{t+1})^\gamma
- \hat f_w (\xi_t)^\gamma \mid \xi_t = \bx]
& \geq \gamma f_w (\bx)^{\gamma-1} \Exp [ f_w (\xi_{t+1}) - f_w (\xi_t) \mid \xi_t = \bx] \\
& ~~ {}+ C f_w (\bx)^{\gamma-2} \Exp [ (f_w (\xi_{t+1}) - f_w (\xi_t))^2 \mid \xi_t = \bx],\end{align*}
for some $C>0$ and any $\bx \in \W_A (\alpha) \cap \W^K (\alpha)$.
Now from Lemma \ref{fmoms} and the conditions
on $\mu(\bx)$ and $\M(\bx)$ it follows that, for some $C>0$, as $\| \bx\| \to \infty$,
\[ \Exp [ \hat f_w (\xi_{t+1})^\gamma
- \hat f_w (\xi_t)^\gamma \mid \xi_t = \bx]
\geq f_w (\bx)^{\gamma-1}  \left[ C f_w (\bx)^{-1} r^{2w-2} + o(r^{w-2}) \right],\]
for
any $\bx \in \W_A (\alpha) \cap \W^K (\alpha)$. Then the result follows since $f_w(\bx)^{-1} \geq r^{-w}$.
 \qed
 
 \subsection{Key estimate}
 \label{sec:key}

Now we state our key lemma for this section. As mentioned above,
the idea is analogous to that used (in a simpler
setting) for  Lemma 6.2 in \cite{bfmp}.

\begin{lemma}
\label{lem2}
Suppose that (A1) and (A2) hold, and
  that for
  some $\sigma^2 \in (0,\infty)$,
  (\ref{sub2}) holds.
Let $\alpha \in (0,\pi]$ and $w = \pi/(2\alpha)$.
 There exist $A \in (0,\infty)$ and $\eps_1,\eps_2>0$
 such that for all $\bx \in \W (\alpha)$
 with $\| \bx \| > A$,
 \[ \Pr   [   \tau_{\alpha,A}    > \eps_1 \| \bx \|^2  \mid \xi_0 = \bx ] \geq \eps_2 \cos (w\varphi)   .\]
\end{lemma}

Our proof   makes repeated use of
  the processes
   $(Y_t (\bx) )_{t \in \Z^+}$  defined for $\bx \in \Z^2$ by
 \begin{equation}
\label{ydef}
Y_t (\bx) := \| \xi_t - \bx \|.
\end{equation}
First note that the triangle inequality implies that
$| Y_{t+1} (\bx) - Y_t (\bx) | \leq \| \xi_{t+1} - \xi_t \|   \leq b$, a.s., 
by (A2).
The next lemma gives
more
information about the increments of $Y_t(\bx)$. For notational ease, for $\bx \in \Z^2$
and $C \in (1,\infty)$ write
\[   S( \bx;C) := \{\by \in \Z^2 : C^{-1} \| \bx \| \leq \| \by \|
 \leq C \| \bx \| \}; ~
   U( \bx;C) := \{\by \in \Z^2 :  \| \by
  -\bx \| \geq C^{-1} \| \bx \| \} .\]

\begin{lemma}
Suppose that (A1) and (A2) hold, and that for
  some $\sigma^2 \in (0,\infty)$, (\ref{sub2}) holds.
  Then for  any $\bx \in \Z^2$ and any $C \in (1,\infty)$, 
  as $\| \bx \| \to \infty$,
   \begin{align}
   \label{ymom0}
    \sup_{\by \in S(\bx;C)}
   \left| \Exp [ Y_{t+1}( \bx)^2 - Y_t (\bx)^2 \mid \xi_t = \by ] - 2 \sigma^2 \right| & = o(1), \\
   \label{ymom2}
    \sup_{\by \in S(\bx;C) \cap U(\bx;C)}
   \left| \Exp [ Y_{t+1}( \bx)  - Y_t (\bx)  \mid \xi_t = \by ] - \frac{1}{2} \sigma^2 \| \by - \bx \|^{-1} \right| &
   = o( \| \bx \|^{-1}), \\
    \label{ymom1}
    \sup_{\by \in S(\bx;C) \cap U(\bx;C)}
       \left| \Exp [ ( Y_{t+1}( \bx)  - Y_t (\bx) )^2  \mid \xi_t = \by ] -  \sigma^2  \right| & = o(1). \end{align}
\end{lemma}
\proof
Conditional on $\xi_t =  \by \in \Z^2$ we have that 
\begin{equation}
\label{eq1}
 {\cal L} ( Y_{t+1} (\bx)   \mid 
\xi_t = \by ) = {\cal L}
( ( \| \by - \bx \|^2  + \| \theta ( \by) \| ^2 +
2 ( \by-\bx)
 \cdot \theta ( \by) )^{1/2}   ) .\end{equation}
Then  (\ref{eq1}) with (\ref{sub2}) yields
 \begin{align*}
  \Exp [ Y_{t+1}(\bx)^2 - Y_t (\bx) ^2 \mid \xi_t = \by ] & = 
   \Exp [ \| \theta (\by) \|^2 ] +
  2 \Exp [ (\by-\bx) \cdot
 \theta (\by)] \nonumber\\
& = 2 \sigma^2 + o(1) + o ( \| \by - \bx \| \| \by \|^{-1} ) = 2 \sigma^2 + o(1),
  \end{align*}
 for   all $\by$ with $C^{-1} \| \bx \| \leq \| \by \| \leq C \| \bx \|$.
This proves (\ref{ymom0}).
Similarly,  by (\ref{eq1}), 
\begin{align}
 \label{ff1}
 \Exp [ Y_{t+1}( \bx)  - Y_t (\bx)  \mid \xi_t = \by ]
  = 
 Y_t(\bx) \Exp \bigg[ \bigg( 1 + \frac{ \| \theta (\by)\|^2 + 2 (\by -\bx) \cdot \theta (\by)}{\| \by - \bx \|^2}
   \bigg)^{1/2} -1 \bigg] .\end{align}
  Taylor's
   theorem applied to the term in square brackets 
   on the right  of (\ref{ff1})
   yields
   \[ \frac{1}{2} \frac{ \| \theta (\by)\|^2 + 2 (\by -\bx) \cdot \theta (\by)}{\| \by - \bx \|^2}
   - \frac{1}{8} \frac{4 ( (\by - \bx) \cdot \theta(\by)) ^2}{\| \by - \bx \|^4 } + O( \| \bx \|^{-3} ) ,\]
   using (A2),
   provided that $C^{-1}  \| \bx \| \leq \| \by -\bx \|$
   and $C^{-1} \| \bx \| \leq \| \by \| \leq C \| \bx \|$.
   Taking expectations of this last expression
  and using (\ref{sub2}), we obtain
   \begin{align*}
 \frac{1}{2} \| \by - \bx \|^{-2} ( 2 \sigma^2 + o(1) )  
- \frac{1}{2} \| \by - \bx \|^{-2} (\sigma^2 + o(1) ) ,\end{align*}
which with (\ref{ff1}) gives   (\ref{ymom2}).
   Finally
 observe that given $\xi_t = \by$,
\begin{align*} (Y_{t+1}(\bx)- Y_t(\bx))^2 
& = (Y_{t+1}(\bx)^2 - Y_t(\bx)^2 ) - 2 \| \by - \bx \| ( Y_{t+1} (\bx) - Y_t (\bx) ) .\end{align*}
So from (\ref{ymom2}) and (\ref{ymom0})
we obtain (\ref{ymom1}).
This completes the proof. \qed\\
  
\noindent
{\bf Proof of Lemma \ref{lem2}.}
  Suppose that $\xi_0 = \bx \in \W (\alpha)$. 
Fix $\alpha' \in (0,\alpha)$, which we will take
close to $\alpha$. First suppose
that $\bx \in \W(\alpha')$, so that the walk
does not start too
close to the boundary of the wedge $\W(\alpha)$.
Note that the shortest distance
from $\bx \in \W(\alpha)$ to the wedge boundary $\partial \W (\alpha)$
is at least
$\| \bx \| \sin (\alpha - | \varphi |)$,
and that for all
$\bx \in \W (\alpha')$, $\varphi \in (-\alpha',\alpha')$
so  this
distance is at least $\eps_0 \| \bx \|$,
where $\eps_0 := \sin (\alpha-\alpha') > 0$.

Suppose that $\by \in B_{ \eps_0 \|\bx\|/2} (\bx) \subset \W(\alpha)$.
Note that for $\by \in B_{ \eps_0 \|\bx\|/2} (\bx)$ we have
\begin{equation}
\label{bz1}
\| \by - \bx \| \leq (\eps_0/2) \| \bx \| ,
~~ \| \by \| \leq ( 1 + (\eps_0/2)) \| \bx \|,
\textrm{~~and~~} \| \by \| \geq ( 1 - (\eps_0/2)) \| \bx \| .
\end{equation}
It then follows from (\ref{ymom0})
and (\ref{bz1}) that
for $\by \in B_{ \eps_0 \|\bx\|/2} (\bx)$,
\begin{align}
\label{kx2}
 \Exp [ Y_{t+1} (\bx)^2 - Y_t (\bx)^2 \mid 
\xi_t =  \by ] 
 = 2 \sigma^2 + o(1) , \end{align}
 as $\| \bx \| \to \infty$.
 For the rest of this proof, let 
 $ \kappa = \min \{ t \in \Z^+ : \| \xi_t - \bx \| \geq \eps_0 \| \bx \|/2 \}$,
 the first exit time of $\Xi$ from $B_{\eps_0 \| \bx \| /2} (\bx)$.
  It follows from (\ref{kx2})
  that for all $\bx \in \W(\alpha')$ with
  $\|\bx\|$
 large enough,  $Y_{t \wedge \kappa} (\bx)^2$ 
 is a nonnegative
  submartingale 
 with respect to the natural
 filtration for $\Xi$,
 and there exists $C \in (0,\infty)$ such that
 for all $\bx \in \W(\alpha')$
 with $\| \bx \|$ sufficiently large and 
 for all $t \in \Z^+$, $\Exp [ Y_{t \wedge \kappa} (\bx) ^2 \mid \xi_0 = \bx
 ] \leq C t \wedge \kappa \leq C t$. 
 
 Then
Doob's submartingale
inequality implies that
there exists $C \in (0,\infty)$ such that
for any $\bx \in \W(\alpha')$ with $\| \bx \|$ sufficiently
large, any $t \in \Z^+$,
and any $x>0$,
\[ \Pr \left[ \left. \max_{0 \leq s \leq t} Y_{s \wedge \kappa} (\bx)
 ^2 \geq x \;\right| \; \xi_0 = \bx
\right]
\leq C t/x .\]
So in time $t = x/(2C)$, there is probability
at least $1/2$ that $\max_{0 \leq s \leq t} Y_{s \wedge \kappa} (\bx)  \leq x^{1/2}$.
Noting that $Y_\kappa (\bx) \geq \eps_0 \| \bx \|/2$ a.s.,
and taking 
$x = \eps_0^2 \| \bx \|^2 / 9$, we conclude that
\[ \Pr \left[ \left. \max_{0 \leq s \leq \eps_0^2 \| \bx \|^2 / (18C)} \| \xi_{s} - \bx \| 
  \leq \eps_0 \| \bx \| /3 \; \right| \; \xi_0 = \bx \right] \geq 1/2.\]
  The event in the last displayed probability
  implies that $\Xi$ remains in $B_{\eps_0 \| \bx \| /2} (\bx) \subset \W(\alpha)$
  till time
  $\eps_0^2 \| \bx \|^2 / (18C)$.
So, for any $\bx \in \W(\alpha')$ with $\| \bx \|$ sufficiently
large,
  \begin{equation}
  \label{prob0}
   \Pr \left[ \left. \tau_{\alpha,A} \geq \frac{\eps_0^2}{18C} \| \bx \|^2 \;\right| \;
    \xi_0 = \bx \right] \geq 1/2.\end{equation}
This yields the statement in the lemma 
for $\bx \in \W(\alpha')$, for any $\alpha' \in (0,\alpha)$.

Now we need to deal with the case $\bx \in \W(\alpha) \setminus \W(\alpha')$.
We take $\alpha' < \alpha$ but close to $\alpha$, so that $\eps_0 = \sin (\alpha-\alpha')$
is small.
 Suppose that
$\bx \in \W(\alpha) \setminus \W(\alpha')$, and, without
loss of generality, that $\varphi >0$; then $\varphi \in [\alpha',\alpha)$.
Set
\[ R:= R(\alpha;\bx) := \begin{cases}
1 & \textrm{ if } \alpha \geq \pi/2 \\
1 \wedge [ (\tan \alpha ) (\cos (\alpha-\varphi) ) ] & \textrm{ if }
 \alpha \in (0,\pi/2) \end{cases} ,\]
and then define
$c(\bx) := \be_r (\alpha) \| \bx \| \cos (\alpha- \varphi) - R (\alpha; \bx) 
\be_\perp (\alpha) \| \bx \|$.

\begin{figure}
\centering
\includegraphics[width=13cm]{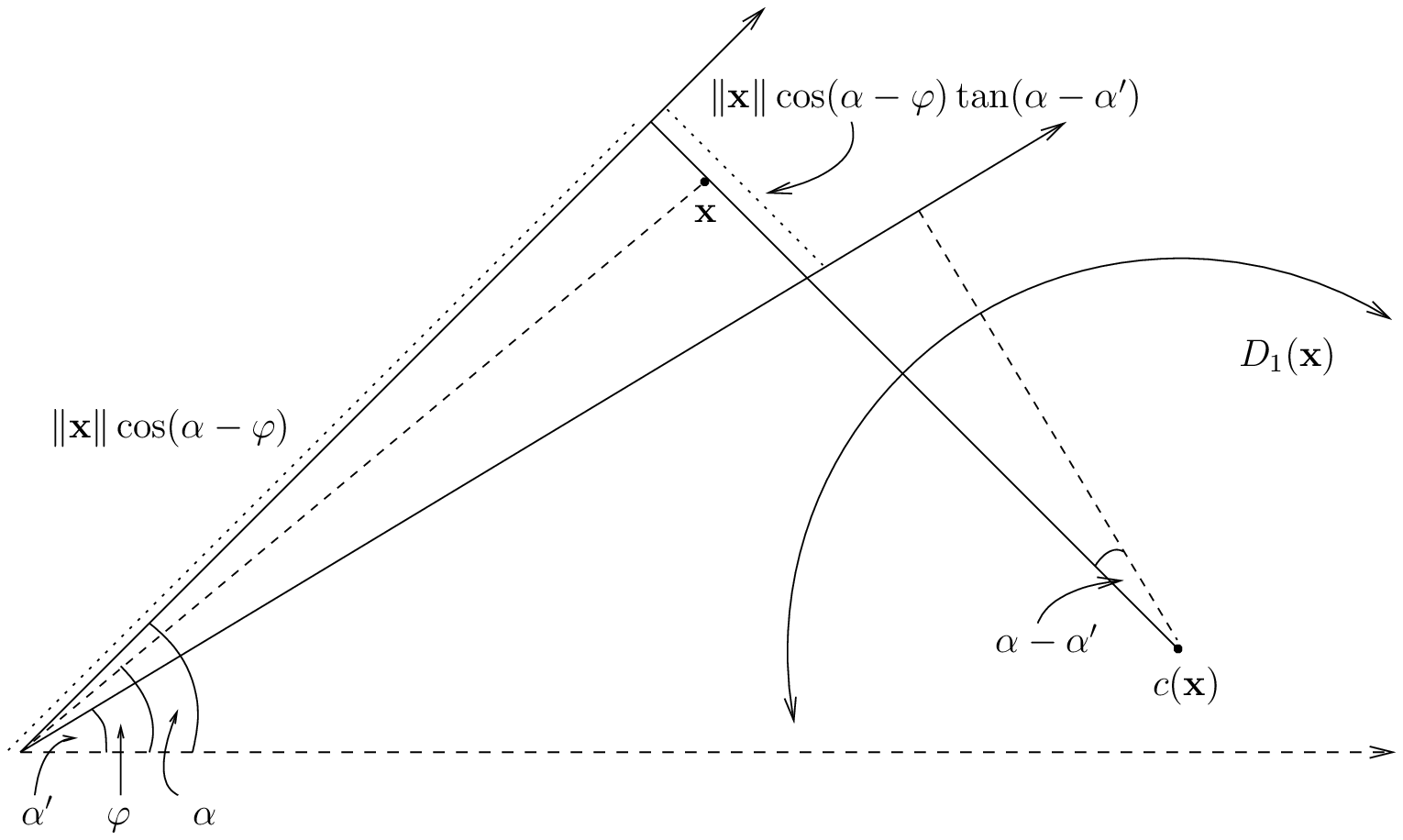}
\caption{The geometrical construction of $c(\bx)$ and $D_1(\bx)$.}
\label{fig3}
\end{figure}

When $R< 1$, this means that $c(\bx) = \be_1 \| \bx \| \cos(\alpha-\varphi) \sec \alpha$
lies on the principal axis of the wedge.
See Figure \ref{fig3} for a typical picture for $R=1$.
Note that
\begin{equation}
\label{cx1}
 \| c(\bx)\| = \| \bx \| \left(R^2 + \cos^2 (\alpha-\varphi) \right)^{1/2} ,\end{equation}
 and also $\bx - c(\bx) = (R-\sin (\alpha-\varphi)) \| \bx \| \be_\perp (\alpha)$, so that
 \begin{equation}
\label{cx2}
 \| \bx - c (\bx) \| = (R- \sin (\alpha -\varphi)) \| \bx\|  \geq \eps_1 \| \bx \|,\end{equation}
 for some $\eps_1>0$ and all $\bx \in \W(\alpha) \setminus \W(\alpha')$ provided
 that $\alpha'$ is close enough to $\alpha$. 
  Also from (\ref{cx1}) we have that for some $\eps_2>0$ and all $\bx \in \W(\alpha) \setminus \W(\alpha')$,
  \begin{equation}
  \label{eq4}
  \eps_2 \| \bx \| \leq \| c (\bx) \| \leq \sqrt{2} \| \bx \| .\end{equation}

Consider the concentric disks
$D_1 (\bx) := B_{R \| \bx \|/2} (c(\bx))$ and $D_2 (\bx) := B_{R \| \bx \|} (c(\bx))$.

If $R=1$, the shortest distance
of $c (\bx)$ from the ray from $\0$
in the $\be_r (\alpha')$
direction is 
\[ \| \bx \| \cos (\alpha - \alpha') - \| \bx \| \sin (\alpha -\alpha') \cos (\alpha - \varphi)
\geq (1- \eps_0) \| \bx \| \cos (\alpha - \alpha') ,\]
for all $\bx \in \W (\alpha) \setminus \W (\alpha')$.
If $R<1$, the corresponding distance is 
$ \| \bx \| \cos (\alpha -\varphi) \sec \alpha \sin \alpha'$. 
In either case,
choosing $\alpha'$ close enough to $\alpha$, it follows
that 
$D_1 (\bx) \subset \W(\alpha')$ for all $\bx \in \W(\alpha) \setminus \W(\alpha')$.
Moreover, for $\eps_0$ small enough,
 for any $\by \in D_2 (\bx)$, by (\ref{cx1}),
\begin{equation}
\label{eq2}
 \| \by \| \geq \| c(\bx) \| - R  \| \bx \| \geq 
\left( (R^2+1-\eps_0^2)^{1/2} - R \right) \| \bx \| \geq \eps_0 \| \bx \|.\end{equation}

We now aim to show that there exists
 $\eps' >0$ such that for
 all $\bx \in \W(\alpha) \setminus \W(\alpha')$
 with $\| \bx\|$ sufficiently large,
\begin{equation}
\label{prob2}
p(\bx) := \Pr \left[ \Xi ~ {\rm visits}~ D_1 (\bx) ~{\rm before~} \R^2 \setminus D_2 (\bx) \mid \xi_0 =\bx \right]
\geq \eps' \cos (w \varphi) .\end{equation}

From the  
 geometrical argument leading up to (\ref{eq2}),
  and the jumps
bound (A2),
it follows that if the event in (\ref{prob2})
occurs,  
$\Xi$ visits a region 
of $\W (\alpha')$ at distance at least
$\eps_0 \| \bx\|$ from $\0$ before
leaving $\W (\alpha)$.
Hence given (\ref{prob2}), (\ref{prob0})
 yields the statement in the lemma in this case also. 
 
 Thus it
 remains to prove (\ref{prob2}). 
With the notation defined at (\ref{ydef}),
we now consider $Y_t (c(\bx)) = \| \xi_t - c(\bx) \|$
 for $\xi_t$
in the annulus
$D_2 (\bx) \setminus D_1 (\bx)$.
For any $\by \in D_2 (\bx) \setminus D_1 (\bx)$ we have
$R \| \bx \|/2 \leq \| \by - c(\bx) \| \leq R \| \bx\|$,
so that
$\| \by \| \leq \| \bx \| + \| c (\bx) \|$. This
together with  (\ref{eq4}) and (\ref{eq2})
implies that 
for $\alpha'$ close enough to $\alpha$ there exist
$C_1, C_2 \in (0,\infty)$ such that
for any $\bx \in \W(\alpha) \setminus \W(\alpha')$
and any $\by \in D_2 (\bx) \setminus D_1 (\bx)$,
\begin{equation}
\label{eq3}
 C_1 \| \bx \| \geq \| \by \| \geq  C_2 \| \bx \|, ~\textrm{and}~ \| \by - c(\bx) \| \geq C_2 \| \bx\|.\end{equation}
Hence by (\ref{eq3}) and (\ref{eq4}),
the estimates (\ref{ymom2}) and
(\ref{ymom1}) are valid
for $Y_t(c (\bx))$ and $\by \in D_2 (\bx) \setminus D_1 (\bx)$,
as $\| \bx \| \to \infty$. Thus we have that there exists $\delta>0$ such
that for $\bx \in \W(\alpha) \setminus \W(\alpha')$
with $\| \bx\|$ large enough and all $\by \in D_2 (\bx) \setminus D_1 (\bx)$,
\begin{align}
\label{ymom1a}
\Exp [ Y_{t+1} (c(\bx)) - Y_t (c(\bx)) \mid \xi_t = \by ]
& = O ( \| \bx \|^{-1} ), \\
\label{ymom2a}
\Exp [ ( Y_{t+1} (c(\bx)) - Y_t (c(\bx)) )^2 \mid \xi_t = \by ]
& > \delta > 0 .\end{align}
 
For $C \in (0,\infty)$ consider now the
process $(Z_t)_{t \in \Z^+}$
defined for $t \in \Z^+$ by 
\[ Z_t := \exp \left\{ C \left( R(\alpha;\bx) - \frac{Y_t(c(\bx))}{\|\bx\|} \right)
\right\} ;\]
then by (\ref{cx2}), $Z_0 = \exp \{ C \sin (\alpha - \varphi) \}$.
Then we have for $t \in \Z^+$ and $\by \in \Z^2$,
\begin{align*}
& \Exp [ Z_{t+1} - Z_t \mid \xi_t = \by ] \\
& =   \exp \left\{ C \left( R-\frac{\|\by - c(\bx)   \|}{\| \bx\|} \right) \right\}
\Exp \left[ \left. \exp \left\{ - \frac{C}{\| \bx \|} ( Y_{t+1}(c(\bx))-Y_t(c(\bx))) \right\}
-1 \; \right| \; \xi_t = \by \right]
.\end{align*}
Since there exist positive
constants $C_1, C_2$ such that $\re^{-x} - 1 \geq - x + C_1 x^2$ for 
 all $x$ with $|x| < C_2$, using the fact that $Y_t(c(\bx))$ has uniformly bounded
 increments we obtain
that for any $\bx \in \W(\alpha) \setminus \W(\alpha')$ and any
$\by \in D_2 (\bx) \setminus D_1 (\bx)$,
\begin{align*} & \Exp \left[ \left. \exp \left\{ - \frac{C}{\| \bx \|} ( Y_{t+1}(c(\bx))-Y_t(c(\bx))) \right\}
-1 \; \right| \; \xi_t = \by \right] \\
& \geq \frac{C}{\| \bx \|} \Exp \left[ \left.
- (Y_{t+1}(c(\bx))-Y_t(c(\bx))) + C_1 \frac{C}{\| \bx \|} (Y_{t+1}(c(\bx))
-Y_t(c(\bx)))^2 \; \right| \; \xi_t = \by \right] .\end{align*}
So by
(\ref{ymom1a}), (\ref{ymom2a})
we may take $C$ large enough such that
for $\xi_0 = \bx \in \W(\alpha) \setminus \W(\alpha')$,
\begin{equation}
\label{ddd} \Exp [ Z_{t+1} - Z_t \mid \xi_t = \by ] \geq 0 ,
\end{equation}
for all $\by \in D_2 (\bx) \setminus D_1 (\bx)$
with $\| \bx \|$ large enough, and all $t \in \Z^+$.
 
 Now to estimate $p(\bx)$ as in (\ref{prob2}),
 we make the
 sets $D_1 (\bx)$ and $\R^2 \setminus D_2 (\bx)$ absorbing.
 Then (using (A2)) $Z_t$ is bounded for this modified 
 random walk, and (using (A1)) $\Xi$ leaves
 $D_2 (\bx) \setminus D_1 (\bx)$ in almost surely finite time.
 Thus as $t \to \infty$, $Z_t$ converges almost surely and
 in $L^1$  to some limit $Z_\infty$ and
 \[ \Exp [ Z_\infty \mid \xi_0 = \bx]
 \leq p(\bx) \exp \{ C R/2 \}   + (1-p(\bx)), \]
 while by (\ref{ddd}) we also have that
 $\Exp [ Z_\infty \mid \xi_0 = \bx]
 \geq \Exp [ Z_0 ] =\exp \{ C \sin (\alpha-\varphi) \}$. 
 Hence  there exists $C \in (0,\infty)$
 such that for all $\bx \in \W(\alpha)
 \setminus \W(\alpha')$ with $\| \bx \|$ large enough
 \[ p(\bx) \geq \frac{ \exp \{ C \sin (\alpha-\varphi) \} -1}
 {\exp \{ C R/2 \} - 1 } \geq \frac{ C}{\re^{CR/2}-1} \sin (\alpha - \varphi) .\]
 Now for $\bx \in \W(\alpha) \setminus \W(\alpha')$ we have that
 $\alpha - \varphi < \alpha -\alpha'$,
 where $\alpha-\alpha'$ is small.
Since, for $a>0$,
 $\frac{\sin (ax)}{\sin (x)} \to a$ as $x \to 0$,
 it follows that there
 exists some $\eps'>0$ such that
 \[ \frac{ C}{\re^{CR/2}-1} \sin (\alpha - \varphi) 
 \geq \eps' \sin (w(\alpha-\varphi)) = \eps' \cos (w \varphi) .\]
 This proves (\ref{prob2}),  
 and so the proof of the lemma is complete.
\qed
  
  \subsection{Proof of Theorem \ref{thm5}(ii)}
  \label{sec:prfnonex}

Now we are ready to complete the proof of Theorem \ref{thm5}(ii). \\
 
\noindent
{\bf Proof of Theorem \ref{thm5}(ii).} Let $\alpha \in (0,\pi]$
and $w = \pi/(2\alpha)$. 
We first
 show that for  $A$ sufficiently large, any $\eps>0$,  and
any $\bx \in \W_A (\alpha)$,
$\Exp [ \tau_\alpha^{(w/2)+\eps} \mid \xi_0 = \bx] =\infty$.
 We proceed in a similar
 way to the proof of Theorem 6.1 in \cite{bfmp}.
 
 Let $A \in (0,\infty)$, to be fixed later.
 For the duration of this proof, to ease notation, set
 $\tau:=\tau_{\alpha,A}$.
 Let $\bx \in \W_A (\alpha)$ be such that
 $f_w (\bx) > A^w$. 
  Suppose, for the purpose
  of deriving a  contradiction, that for some $\eps>0$,
$\Exp [  \tau^{(w/2)+\eps} \mid \xi_0 = \bx] < \infty$.
Let $\Xi' = (\xi'_t)_{t \in \Z^+}$ be an independent
copy of $\Xi$, and let $\tau'$ be the corresponding
independent copy of $\tau$.
Then for any $t \in \N$, by conditioning on
$\xi_t$ and using the Markov property at time $\tau$, 
\begin{align*} \Exp [   \tau ^{(w/2) +\eps} \mid \xi_0 = \bx ]   \geq 
  \Exp \left[ \Exp 
\left[ ( t+  \tau' )^{(w/2) +\eps} 
 \mid \xi'_0 = \xi_t \right] \1_{\{  \tau > t\}}  \mid \xi_0 = \bx \right].\end{align*}
Hence by Lemma \ref{lem2}, for $A$ large enough,
\begin{align*} \Exp [  \tau ^{(w/2) +\eps} \mid \xi_0 =\bx] &
\geq  \eps_2  \Exp [ (t+\eps_1 \| \xi_t \| ^2  )^{(w/2) +\eps}  \cos  (w\varphi(\xi_t))
  \1_{\{  \tau  > t\}}  \mid \xi_0 =\bx] \\
& \geq   
C \Exp
\left[ (\hat f_w ( \xi_{t \wedge \tau} ) )^{1 +(2/w) \eps}   \mid \xi_0 = \bx
\right] - A^{w+2\eps},
\end{align*}
for some $C \in (0,1)$, any $t \in \N$,
using the fact that $\hat f_w (\xi_{\tau})   \leq A^w$ a.s..
Thus
 under the hypothesis
$\Exp [  \tau^{(w/2)+\eps} \mid \xi_0 = \bx] < \infty$,
 for some $\eps'>0$
the process $(\hat f_w (\xi_{t \wedge   \tau } ) )^{1+\eps'}$
 is uniformly integrable. It follows
 (since by hypothesis 
  $\tau  < \infty$ a.s.) 
 that 
as $t \to \infty$,
$\Exp [ ( \hat f_w (\xi_{t \wedge   \tau } ) )^{1+\eps'} \mid \xi_0 =\bx ] 
\to \Exp [  ( \hat f_w (\xi_{  \tau } )  )^{1+\eps'} \mid \xi_0 =\bx] \leq A^{w(1+\eps')}$. However, by
the submartingale property
(\ref{sub3}),   for $A$ large enough,
$\Exp [ ( \hat f_w (\xi_{t \wedge \tau  } ) )^{1+\eps'} \mid \xi_0 =\bx   ] \geq    ( \hat f_w (\bx )   )^{1+\eps'}
> A^{w(1+\eps')}$ for all $t \in \N$, given our condition on $\bx$.  
Thus we have the desired contradiction,
and $\Exp [  \tau^{(w/2)+\eps} \mid \xi_0 = \bx] = \infty$
for any $\bx \in \W_A(\alpha)$ with $f_w (\bx) > A^w$. 
Since $\tau_\alpha \geq \tau$ a.s., this implies that 
$\Exp [  \tau_\alpha^{(w/2)+\eps} \mid \xi_0 = \bx] = \infty$
for any $\bx \in \W_A(\alpha)$ with $f_w (\bx) > A^w$.  
  Lemma \ref{lem88} extends the conclusion
to any $\bx \in \W(\alpha)$ with $\| \bx \|$ large enough.
\qed

\section*{Acknowledgements}

Some of this work was done while AW was at the University of Bristol,
supported by the Heilbronn Institute for Mathematical Research.


\begin{thebibliography}{99}
 
\bibitem{aim} S. Aspandiiarov, R. Iasnogorodski, and M. Menshikov,
Passage-time moments for nonnegative stochastic processes and an application to reflected random walks in a quadrant,
{\em Ann. Probab.} {\bf 24} (1996) 932--960. 
   
\bibitem{bfmp} V. Belitsky, P.A. Ferrari, M.V. Menshikov, and S.Yu. Popov, A mixture
of the exclusion process and the voter model, 
{\em Bernoulli} {\bf 7} (2001) 119--144.

\bibitem{burkh} D.L. Burkholder, Exit times of Brownian motion,
harmonic majorization, and Hardy spaces, 
{\em Adv. Math.} {\bf 26} (1977) 182--205.
  
 \bibitem{cohen} J.W. Cohen, On the random walk with zero drifts in
 the first quadrant of $\R_2$,
 {\em Comm. Statist. Stochastic Models} {\bf 8} (1992) 359--374.
 
 \bibitem{dante} R.D. DeBlassie, Exit times from cones in $\R^n$ 
 of Brownian motion,
 {\em Probab. Theory Relat. Fields} {\bf 74} (1987) 1--29.
   
 \bibitem{fmm} G. Fayolle, V.A. Malyshev, and M.V. Menshikov,
Topics in the Constructive Theory of Countable Markov Chains,
Cambridge University Press, 1995.

\bibitem{foster}
F.G. Foster, On the stochastic matrices associated with certain queuing processes,
{\em Ann. Math. Statist.} {\bf 24} (1953) 355--360.

\bibitem{flp} 
J.-D. Fouks, E. Lesigne, and M. Peign\'e,
\'Etude asymptotique d'une marche al\'eatoire centrifuge,
{\em Ann. Inst. H. Poincar\'e Probab. Statist.} {\bf 42} (2006) 147--170.

\bibitem{fukai} Y. Fukai,
Hitting time of a half-line by two-dimensional
random walk,
{\em Probab. Theory Relat. Fields} {\bf 128} (2004)
323--346.

\bibitem{garbit} R. Garbit, Temps de sortie d'une c\^one pour une marche al\'eatoire centr\'ee,
{\em C. R. Math. Acad. Sci. Paris} {\bf 345} (2007) 587--591.

\bibitem{kingman} J.F.C. Kingman,
The ergodic behaviour of random walks,
{\em Biometrika} {\bf 48} (1961) 391--396.

\bibitem{klein} L.A. Klein Haneveld and A.O. Pittenger,
Escape time for a random walk from an orthant,
{\em Stochastic Processes Appl.} {\bf 35} (1990) 1--9.
 
\bibitem{lamp1} J. Lamperti, 
Criteria for the recurrence and
transience of stochastic processes I, 
{\em J. Math. Anal. Appl.} {\bf 1} (1960) 314--330.
   
\bibitem{lamp2} J. Lamperti, Criteria for
stochastic processes II: passage-time moments,
{\em J. Math. Anal. Appl.} {\bf 7} (1963) 127--145.

 \bibitem{lawler91} G.F. Lawler, Estimates for differences and Harnack
 inequality for difference operators coming from random walks with
 symmetric, spatially inhomogeneous, increments,
 {\em Proc. London Math. Soc.} {\bf 63} (1991) 552--568.
   
   \bibitem{lawlerbook} G.F. Lawler,
   Intersections of Random Walks,
   {\em
   Probability and Its Applications}, Birkh\"auser, Boston,
   1996.
   
  \bibitem{mmw1} I.M. MacPhee, M.V. Menshikov, and A.R. Wade,
Angular asymptotics for multi-dimensional non-homogeneous random walks with asymptotically zero drift,
{\em Markov Process. Relat. Fields} {\bf 16} (2010) 351--388.
   
   \bibitem{malyshev} V.A. Malyshev, Classification
   of two-dimensional positive
   random walks and almost
   linear semimartingales,
  {\em Soviet Math. Dokl.} {\bf 13} (1972) 136--139;
  translated from {\em Dokl. Akad. Nauk SSSR} {\bf 202} (1972) 526--528 (in Russian).
       
  \bibitem{mvw} M.V. Menshikov, M. Vachkovskaia, and A.R. Wade,
  Asymptotic behaviour of randomly reflecting billiards
  in unbounded tubular domains,
 {\em J. Statist. Phys.} {\bf 132} (2008) 1097--1133.
    
    \bibitem{mw3} M.V. Menshikov and A.R. Wade,
    Rate of escape and central limit theorem for the supercritical Lamperti
    problem, {\em Stochastic Process. Appl.} {\bf 120} (2010) 2078--2099. 
    
  \bibitem{mustapha} S. Mustapha,
  Gaussian estimates for spatially inhomogeneous random walks
  on $Z^d$,
  {\em Ann. Probab.} {\bf 34} (2006) 264--283.
     
  \bibitem{spitzer} F. Spitzer, Some theorems concerning $2$-dimensional
  Brownian motion,
  {\em Trans. Amer. Math. Soc.} {\bf 87} (1958) 187--197.
  
  \bibitem{spitzerrw} F. Spitzer, Principles of Random Walk,
  2nd edition, Springer, New York, 1976.
  
\bibitem{var1} N.Th. Varopoulos, 
Potential theory in conical domains, 
{\em Math. Proc. Cambridge Philos. Soc.} {\bf 125} (1999) 335--384. 

\bibitem{var2} N.Th. Varopoulos, 
Potential theory in conical domains II,
{\em Math. Proc. Cambridge Philos. Soc.} {\bf 129} (2000) 301--319.
  
\end{thebibliography}
\end{document}